\documentclass[preprint,12pt]{elsarticle}

\usepackage{amssymb,amsmath,amsfonts,pbox,subfig,graphicx,mwe}

  \def\clap#1{\hbox to 0pt{\hss#1\hss}}

\usepackage{bm}
\providecommand{\mat}[1]{\bm{#1}}%
\renewcommand{\vec}[1]{\mathbf{#1}}

\providecommand{\mA}{\ensuremath{\mat{A}}}
\providecommand{\mB}{\ensuremath{\mat{B}}}
\providecommand{\mC}{\ensuremath{\mat{C}}}

\providecommand{\mW}{\ensuremath{\mat{W}}}

\providecommand{\vb}{\ensuremath{\vec{b}}}

\providecommand{\vf}{\ensuremath{\vec{f}}}

\providecommand{\vw}{\ensuremath{\vec{w}}}
\providecommand{\vx}{\ensuremath{\vec{x}}}

\providecommand{\vnu}{\bm{\nu}}

\newcommand{\hmC}{\hat{\mC}}

\newcommand{\hLambda}{\hat{\Lambda}}

\newcommand{\hmW}{\hat{\mW}}

\newcommand{\sS}{\mathcal{S}}

\newcommand{\Exp}[1]{\mathbb{E}\left[#1\right]}

\newcommand{\Var}[1]{\operatorname{Var}\left[#1\right]}

\newcommand{\bmat}[1]{\begin{bmatrix}#1\end{bmatrix}}

\newcommand{\diag}[1]{\operatorname{diag}(#1)}

\usepackage{url}
\usepackage{algorithm}
\usepackage{color}
\newtheorem{theorem}{Theorem}[section]

\newproof{pf}{Proof}

\journal{Reliability Engineering \& System Safety}

\begin{document}

\begin{frontmatter}



\title{Global sensitivity metrics from active subspaces}


\author{Paul G.~Constantine\corref{cor1}\fnref{CSM}}
\cortext[cor1]{Corresponding author}
\address[CSM]{Department of Applied Mathematics and Statistics, Colorado School of Mines, Golden, CO}
\ead{paul.constantine@mines.edu}
\ead[url]{http://inside.mines.edu/~pconstan}

\author{Paul Diaz\fnref{CU}}
\address[CU]{Department of Aerospace Engineering, University of Colorado Boulder, Boulder, CO}

\begin{abstract}
Predictions from science and engineering models depend on several input parameters. Global sensitivity analysis quantifies the importance of each input parameter, which can lead to insight into the model and reduced computational cost; commonly used sensitivity metrics include Sobol' total sensitivity indices and derivative-based global sensitivity measures. \emph{Active subspaces} are an emerging set of tools for identifying important directions in a model's input parameter space; these directions can be exploited to reduce the model's dimension enabling otherwise infeasible parameter studies. In this paper, we develop global sensitivity metrics called \emph{activity scores} from the active subspace, which yield insight into the important model parameters. We mathematically relate the activity scores to established sensitivity metrics, and we discuss computational methods to estimate the activity scores. We show two numerical examples with algebraic functions taken from simplified engineering models. For each model, we analyze the active subspace and discuss how to exploit the low-dimensional structure. We then show that input rankings produced by the activity scores are consistent with rankings produced by the standard metrics. 
\end{abstract}


\end{frontmatter}


\section{Introduction}
\label{sec:intro}

\noindent Complex science and engineering models depend on several input parameters. The complexity often creates uncertainty about which parameters are most important for describing physical processes and model predictions. Sensitivity analysis seeks to identify the most important parameters for particular outputs of interest; effectively identifying important parameters leads to insight into the model~\cite{saltelli2008global}. Sensitivity analysis may enable the scientist to ignore less important parameters, thus reducing the cost of parameter studies---e.g., optimization or uncertainty quantification---that are essential for confident predictions. Sensitivity metrics are broadly classified as \emph{local}, which measure the model's response to small perturbations around a nominal parameter value, and \emph{global}, which assess the importance of each variable over a range of parameters; we are primarily concerned with global sensitivity metrics. The \emph{total sensitivity indices} of Sobol' are some of the most widely used metrics for global sensitivity analysis~\cite{sobol2001global}. Other common metrics average local gradients to build global metrics~\cite{kucherenko2009monte}. The appropriate choice of metric depends on the application. 

Using sensitivity metrics to ignore less important parameters is a form of dimension reduction in the model's input parameter space. This type of reduction is limited to the coordinates induced by the model parameters. Our recent work in \emph{active subspaces} generalizes such parameter space dimension reduction to directions in the space of model inputs~\cite{Constantine2014,asm2015}. A model's output might depend strongly on each coordinate but only through a limited number of directions. Each direction corresponds to a linear combination of the input parameters that we call an \emph{active variable}. Our work has developed methods for identifying an active subspace~\cite{constantine2014computing} and exploiting it, if present, to reduce the dimension of subsequent parameter studies~\cite{Constantine2014}. We have applied this approach to aerospace shape optimization~\cite{Lukaczyk2014}, automotive design~\cite{Othmer2016}, uncertainty quantification for multiphysics scramjet models~\cite{constantine2014exploiting}, sensitivity analysis in integrated hydrologic models~\cite{Jefferson2015}, global spatial sensitivity in subsurface flow~\cite{Gilbert2016}, dimension reduction in single-diode solar cell models~\cite{constantine2015discovering}, uncertainty quantification for combustion models~\cite{Bauerheim2016}, time-dependent sensitivity analysis for lithium ion battery models~\cite{constantine2016battery}, and sensitivity analysis for Ebola virus population dynamics~\cite{diaz2016modified}. The number of input parameters in these models varies between 5 and 100. The active subspace quantifies low-dimensional structure in the input/output map of a given model. The aforementioned list of science and engineering applications demonstrates that such structures are present in many real world models; moreover we have shown that these low-dimensional structures are generically present in any physical system with units~\cite{constantine2016many}. If an active subspace is sufficiently low dimensional in a given model, then the researcher may exploit the active subspace to enable otherwise infeasible parameter studies. The way in which she exploits the structure depends on the study. For example, to find the range of possible outputs over the parameter space (i.e., minimization and maximization), she may exploit a global, monotonic trend in the relationship between the active variables and the output quantity of interest; we apply this approach in~\cite{constantine2014exploiting} to find the safe operating in multiphysics scramjet model. Alternatively, if she wishes to construct a cheap response surface that mimics the model's input/output relationship, she may train the surface on only the active subspace coordinates; see~\cite{constantine2016near} for details. 

One difficulty in subspace-based dimension reduction is how to draw insight into the original model parameters once an active subspace has been identified. For example, how might a climate scientist use the fact that $3\times(\text{cloud parameter}) + 2\times(\text{radiation parameter})$ is an important linear combination in a climate model? In this paper, we propose global sensitivity metrics based on active subspaces; our main goal is to demonstrate that the proposed metrics function like standard sensitivity metrics---even though the proposed metrics' derivation and motivation are quite different. Toward this goal, we mathematically relate the proposed metrics to existing sensitivity metrics: Sobol' total sensitivity indices and derivative-based global sensitivity measures. We also discuss two numerical procedures for estimating the proposed metrics: tensor product Gauss quadrature and Monte Carlo. For the Monte Carlo method, we propose bootstrap-based standard error estimates. We use two simplified engineering models to compare the proposed metrics to existing metrics: (i) a 7-parameter algebraic model of piston cycle time and (ii) a 6-parameter algebraic model of a transformerless circuit voltage. The computationally inexpensive models with relatively few parameters allow us to perform a thorough numerical error study, where we compare the Monte Carlo-based estimates to quadrature-based reference values. We also numerically compare the proposed metrics to the coefficients of a least-squares fit linear model; the related standard regression coefficients are common sensitivity metrics. We show that the rankings produced by all metrics are similar for both models; this result supports the claim that the active subspace-based sensitivity metrics are comparable to the standard metrics. It is not our goal to demonstrate that the active subspace-based metrics are superior to existing metrics---either in characterizing sensitivity for a given model or in computational cost for estimating the metrics. In fact, when derivatives are available, the derivative-based global sensitivity measures use a comparable number of gradient evaluations, and they do not require an eigendecomposition.

The remainder of the paper proceeds as follows. Section \ref{sec:background} reviews three common sensitivity metrics: Sobol' total sensitivity indices, derivative-based global sensitivity measures, and linear model coefficients. This section sets up notation for the remainder of the paper. Section \ref{sec:as} reviews active subspaces. Section \ref{sec:safromas} mathematically defines the proposed sensitivity metrics derived from the active subspace and relates the proposed metrics to the total sensitivity indices and the derivative-based global sensitivity measures. Section \ref{sec:exp} shows the results of the numerical experiments. And Section \ref{sec:conc} summarizes the results, discusses limitations of the proposed metrics, and offers directions for future research.

\subsection{Setup and notation}
\label{sec:setup}

\noindent To study the sensitivity metrics derived from active subspaces, we consider a differentiable and square-integrable function $f=f(\vx)$ defined on an $m$-dimensional hypercube $\vx=[x_1,\dots,x_m]^T\in[-1,1]^m$. Define the uniform probability density function $\rho(\vx)=2^{-m}$ for $\vx\in[-1,1]^m$ and zero elsewhere. The mean and variance of $f$ are
\begin{equation}
\Exp{f} \;=\; \int f(\vx)\,\rho(\vx)\,d\vx,\qquad 
\Var{f} \;=\; \int \big(f(\vx)-\Exp{f}\big)^2\,\rho(\vx)\,d\vx.
\end{equation}
We emphasize that $f$ is a deterministic function; its mean and variance do not represent any probabilistic information. Denote the gradient of $f$ by $\nabla f(\vx)\in\mathbb{R}^m$ with partial derivatives $\frac{\partial f}{\partial x_i}(\vx)$, and assume the partial derivatives are square-integrable over the hypercube domain. All vector norms are the standard Euclidean norm, and matrix norms are the matrix 2-norm (i.e., the spectral norm).

This set up is an idealized input/output map. To formally apply the following methods and results to a real application, the researcher must (i) identify a scalar-valued quantity of interest that satisfies the assumptions (e.g., square-integrable and differentiable), (ii) identify a vector of continuous input parameters, and (iii) transform the input parameter space to the hypercube with the uniform measure. 

\section{Global sensitivity metrics}
\label{sec:background}

\noindent We briefly review three global sensitivity metrics used in practice: variance-based total sensitivity indices, derivative-based global sensitivity measures, and least-squares-fit linear model coefficients. In each case, we (i) define the metric, (ii) offer a method for high accuracy computation with numerical quadrature, (iii) detail a Monte Carlo method for low accuracy estimation, and (iv) mention practical error estimates for the Monte Carlo estimates. We use lower case Greek letters to denote the various sensitivity metrics, where each symbol's subscript identifies its corresponding variable (i.e., component of the vector of inputs $\vx$). The book by Saltelli, et al.~\cite{saltelli2008global} provides an overview of global sensitivity analysis and its metrics including common metrics we do not review, such as correlation coefficients.

\subsection{Sobol' total sensitivity indices}
\label{sec:variance}

\noindent The total sensitivity indices were introduced by Sobol'~\cite{sobol2001global}; several papers describe these metrics, propose extensions~\cite{Owen13}, and develop efficient computational procedures~\cite{Liu2006,Sobol2008,Sudret2008,crestaux2009polynomial,Saltelli2010,Weirs2012}. The total sensitivity indices are derived from a particular decomposition of the square-integrable $f(\vx)$ called the functional ANOVA decomposition or the variance-based decomposition; we use the concise notation from Liu and Owen~\cite{Liu2006} to describe the decomposition and the derived sensitivity indices. Let $u$ be a subset of the integers $\{1,\dots,m\}$ with cardinality $|u|$. Define $\vx^u\in[-1,1]^{|u|}$. The subset $u_c$ is the complement of $u$ in $\{1,\dots,m\}$ with cardinality $|u_c| = m-|u|$. The ANOVA decomposition can be written
\begin{equation}
\label{eq:anova}
f(\vx) \;=\; \sum_{u\subseteq\{1,\dots,m\}} f_u(\vx),
\end{equation}
where $f_u$ depends on $\vx$ only through $\vx^u$. Each function $f_u$ is defined as
\begin{equation}
f_u(\vx) \;=\; 2^{-|u_c|}\,\int_{[-1,1]^{|u_c|}} f(\vx)\,d\vx^{u_c} - \sum_{v\subsetneq u} f_v(\vx),
\end{equation}
where $f_{\emptyset} = \Exp{f}$. In the standard construction, the functions $f_u$ are mutually orthogonal. The variance of $f$ can be decomposed as
\begin{equation}
\label{eq:varfu}
\Var{f} \;=\; \sum_{u\subseteq\{1,\dots,m\}} \Var{f_u}. 
\end{equation}
Define the set $\sS_i$ to be the set of subsets of $\{1,\dots,m\}$ containing the index $i$. The total sensitivity index $\tau_i$ for the $i$th variable $x_i$ can be written as sums of subsets of the partial variances $\Var{f_u}$ divided by the total variance,
\begin{equation}
\label{eq:tsi}
\tau_i \;=\; \frac{\sum_{u\subseteq\sS_i} \Var{f_u}}{\Var{f}}.
\end{equation}
Some authors use the notation $S_{T_i}$~\cite{saltelli2008global} or $S^{\text{tot}}_i$~\cite{sobol2001global} to denote the total sensitivity indices. We use the lower case Greek letter $\tau_i$ for consistency. Other sensitivity metrics can be defined from the variances $\Var{f_u}$; we focus exclusively on $\tau_i$ from \eqref{eq:tsi} to compare with the metrics derived from active subspaces; see Section \ref{sec:safromas}.

When the dimension of $\vx$ is sufficiently low, $f(\vx)$ is sufficiently smooth, and $f(\vx)$ can be evaluated quickly, one can accurately estimate $\tau_i$ using a high order numerical quadrature rule, such as a tensor product Gauss-Legendre quadrature rule. This approach is equivalent to using the Legendre series' coefficients as in~\cite{Sudret2008,crestaux2009polynomial}, which we use to compute reference values for the numerical examples in Section \ref{sec:exp}. 

If $f(\vx)$ does not satisfy these assumptions (e.g., if the dimension of $\vx$ is too large to permit tensor product quadrature), one can use a less accurate Monte Carlo method as in~\cite[Section 4.6]{saltelli2008global}. Several alternative Monte Carlo methods have been proposed; see~\cite{Saltelli2010,Weirs2012} for a detailed discussion. We review Jansen's method~\cite{Jansen1999}, which is the method we use for numerical experiments, and a bootstrap standard error. First draw $2M$ independent samples $\vx_i$ according to the uniform density $\rho(\vx)$ on $[-1,1]^m$. Arrange these samples into two $M\times m$ matrices, $\mA$ and $\mB$, defined as
\begin{equation}
\mA = \bmat{\vx^T_1\\ \vdots\\ \vx^T_M},\qquad \mB = \bmat{\vx^T_{M+1}\\ \vdots\\ \vx^T_{2M}}.
\end{equation}
Define the $M\times m$ matrix $\mB_i$ to be the matrix $\mB$ with the $i$th column replaced by the $i$th column of $\mA$. Define the $M$-vector $\vf_{\mA}$ as
\begin{equation}
\vf_{\mA} \;=\; \bmat{f(\vx_1)\\ \vdots \\ f(\vx_M)}.
\end{equation}
In words, the vector $\vf_{\mA}$ contains evaluations of the function $f$ at each row of the matrix $\mA$. Similarly, define the $M$-vectors $\vf_{\mB_i}$ whose elements contain evaluations of $f$ at the rows of their respective subscripted matrices. These calculations use $M(1+m)$ unique evaluations of $f$. The total sensitivity indices may be estimated with Jansen's formula~\cite{Jansen1999,Saltelli2010},
\begin{equation}
\label{eq:saltelli_tsi}
\tau_i \;\approx\; \hat{\tau}_i \;=\;
\frac{\|f_{\mA} - f_{\mB_i}\|^2}{2\,M\,\hat{\sigma}_f^2},
\end{equation}
where $\hat{\sigma}_f^2$ is the sample variance of the $M$ evaluations of $f$ in $\vf_{\mA}$, which approximates the true variance of $f$. 

Under certain assumptions on the partial variances $\Var{f_u}$ from \eqref{eq:varfu}, one can derive the standard error of the Monte Carlo estimator from a chi-squared distribution~\cite{Jansen1999}. However, these assumptions are not satisfied for a general nonlinear $f(\vx)$. To estimate the standard error, we employ a nonparametric bootstrap~\cite{efron1994introduction}. We denote the bootstrap standard error for the Monte Carlo estimate $\hat{\tau}_i$ by $\widehat{se}_{\tau_i}$. And we numerically study the bootstrap standard error as a function of the number of samples in Section \ref{sec:exp}.


\subsection{Derivative-based global sensitivity measures}
\label{sec:derivative}

\noindent A natural notion of sensitivity is how a model output responds to small changes in its inputs, i.e., a derivative. However, derivatives only provide information at a single parameter value. Morris developed sensitivity metrics based on sampling and averaging coarse finite difference approximations to derivatives~\cite{Morris1991}. Kucherenko, et al.~\cite{kucherenko2009monte} extended this idea to continuous derivatives and integration, where the the averages are estimated with Monte Carlo methods. We focus on one particular derivative-based sensitivity metric studied by Sobol' and Kucherenko~\cite{Sobol2009}:
\begin{equation}
\label{eq:avgsqder}
\nu_i \;=\; \int 
\left(\frac{\partial f}{\partial x_i}(\vx)\right)^2
\,\rho(\vx)\,d\vx,\qquad i=1,\dots,m.
\end{equation}
This metric averages local sensitivity information to make it a global metric. We discuss this metric's connection to active subspaces in Section \ref{sec:safromas}. 

If $f$'s derivatives are sufficiently smooth, and if the dimension of $\vx$ is low, then one can accurately estimate $\nu_i$ using high order numerical quadrature, such as tensor-product Gauss-Legendre quadrature; we use the quadrature approach for reference values in the numerical examples in section \ref{sec:exp}.

The Monte Carlo estimate of $\nu_i$ using $M$ samples is
\begin{equation}\label{eq:dbsm_mc}
\nu_i \; \approx \; \hat{\nu}_i = \frac{1}{M}\sum_{j=1}^M \left(\frac{\partial f}{\partial x_i}(\vx_j)\right)^2,
\end{equation}
where the $\vx_j$'s are drawn independently according to $\rho(\vx)$. The standard error $se_{\nu_i}$ is
\begin{equation}
se_{\nu_i} \;=\;  \left[ 
\frac{1}{M-1} \sum_{j=1}^M
\left(
\left(\frac{\partial f}{\partial x_i}(\vx_j)\right)^2 - \hat{\nu}_i
\right)^2 
\right]^{1/2}.
\end{equation}
In section \ref{sec:exp}, we study the standard error of the estimator by investigating $se_{\nu_i}$ for different values of $M$.

\subsection{Linear model coefficients}
\label{sec:regression_coefficients}

\noindent A simple and easily computable sensitivity metric can be derived from the coefficients of a least-squares fit linear approximation to $f(\vx)$. Such a metric is valid when $f$ is smooth and monotonic along each component of $\vx$. Saltelli, et al.~\cite[Section 1.2.5]{saltelli2008global} discuss a similar metric in connection with a linear regression model. Since $f(\vx)$ is deterministic, the statistical interpretation of the regression coefficients is not valid. Nevertheless, we can devise a Monte Carlo-based least-squares procedure to estimate the coefficients of a linear approximation as follows. Draw $\vx_j$ independently according to $\rho(\vx)$ for $j=1,\dots,M$, and compute $f_j=f(\vx_j)$. Next, use least-squares to fit the intercept $b_0$ and coefficients $\vb=[b_1,\dots,b_m]^T$ of a linear approximation,
\begin{equation}
\label{eq:linmod}
f_j \;\approx\; b_0 + \vb^T\vx_j,\qquad j=1,\dots,M.
\end{equation}
Scaling the coefficients $b_i$ as follows produces useful sensitivity metrics called the \emph{standardized regression coefficients}; see~\cite[Section 1.2.5]{saltelli2008global}. Define $\beta_i$ as
\begin{equation}\label{eq:srq}
\beta_i \;=\; \frac{b_i}{\sqrt{3}\,\hat{\sigma}_f},
\end{equation}
where $\hat{\sigma}_f$ is the standard deviation of the set of samples $\{f_j\}$. The $1/\sqrt{3}$ comes from scaling the coefficients by the variance of the uniformly distributed random variable $x_i$ with support $[-1,1]$, which is consistent with the problem set up in Section \ref{sec:setup}.  If $f(\vx)$ is monotonic and sufficiently smooth, then $\beta_i$ gives a rough indication of how $f$ changes in response to a change in $x_i$ over the space of inputs. In contrast to the total sensitivity index $\tau_i$ from \eqref{eq:saltelli_tsi} and the derivative-based measures $\nu_i$ from \eqref{eq:avgsqder}, the $\beta_i$'s are signed, so they indicate if $f$ will increase or decrease given a positive change in $x_i$. 

We can compute reference values for $\beta_i$ by noting the relationship between the $m$-dimensional Legendre series truncated to degree 1 (i.e., only a constant and $m$ linear terms) and the least-squares linear approximation \eqref{eq:linmod}; the Legendre series is often referred to as the \emph{generalized polynomial chaos} associated with the uniform density~\cite{Sudret2008,Xiu2002}. In particular, the degree-1 Legendre approximation is the best linear approximation in the $L_2$ norm. Its coefficients admit a simple integral representation by virtue of the Legendre polynomials' orthogonality. Thus, we can estimate the coefficients with high order Gauss quadrature~\cite{Constantine2012} and scale them by a quadrature-based estimate of the variance to get the coefficients of a linear monomial approximation as in \eqref{eq:linmod}.

When reference values are too expensive to compute, we can estimate the standard error in $\beta_i$ with a bootstrap standard error, which we denote $\widehat{se}_{\beta_i}$. See \cite[Section 9.5]{efron1994introduction} for a detailed description of using the bootstrap for linear regression coefficients. Note that randomness in the estimates is a result of the Monte Carlo method; there is no randomness in the evaluation of $f(\vx)$. In section \ref{sec:exp}, we numerically study the decay in the bootstrap standard error for different values of $M$.

\section{Active subspaces}
\label{sec:as}

\noindent The active subspace~\cite{asm2015} is defined by the eigenvectors of the following $m\times m$ symmetric, positive semidefinite matrix,
\begin{equation}
\label{eq:C}
\mC \;=\; \int \nabla f(\vx)\,\nabla f(\vx)^T\,\rho(\vx)\,d\vx \;=\;
\mW\Lambda\mW^T,
\end{equation}
where $\mW=[\vw_1,\dots,\vw_m]$ is the orthogonal matrix of eigenvectors, and $\Lambda=\diag{\lambda_1,\dots,\lambda_m}$ is the diagonal matrix of eigenvalues in decreasing order. The eigenpairs are functionals, i.e., properties of the given $f$, similar to the total sensitivity indices or derivative-based measures. The eigenpairs satisfy
\begin{equation}
\label{eq:eigenval}
\lambda_i \;=\; \vw_i^T\mC\vw_i \;=\; 
\int \big(\nabla f(\vx)^T\vw_i \big)^2\,\rho(\vx)\,d\vx,
\end{equation}
which implies that $\lambda_i=0$ if and only if $f$ is constant along the direction $\vw_i$. If such structure is present in the given $f$, then the structure can be exploited to reduce the parameter space dimension for parameter studies. Suppose $\lambda_n>\lambda_{n+1}$ for some $n<m$. Then we can partition the eigenpairs,
\begin{equation}
\label{eq:partition}
\Lambda = \bmat{\Lambda_1 & \\ & \Lambda_2},\qquad
\mW = \bmat{\mW_1 & \mW_2},
\end{equation}
where $\Lambda_1$ contains the first $n$ eigenvalues, and $\mW_1$ is the $m\times n$ matrix containing the first $n$ eigenvectors. The \emph{active subspace} of dimension $n$ is the column span of $\mW_1$. Not all square-integrable and differentiable functions admit an active subspace. For example, for the function $f(\vx)=\|\vx\|^2$, which is radially symmetric, and a uniform density $\rho(\vx)$, the eigenvalues are are all equal. 

If $\lambda_{n+1},\dots,\lambda_m$ are sufficiently small, then a reasonable low-dimensional approximation of $f$ is
\begin{equation}
\label{eq:g}
f(\vx) \;\approx\; g(\mW_1^T\vx),
\end{equation}
where $g:\mathbb{R}^n\rightarrow\mathbb{R}$. In~\cite{Constantine2014}, we propose a particular form for $g$ based on a conditional average, and we bound its mean-squared error in terms of the eigenvalues $\lambda_{n+1},\dots,\lambda_m$. 


To identify an active subspace, we must estimate the eigenpairs $\Lambda$, $\mW$ and examine the estimated eigenvalues. If the dimension $m$ is sufficiently small, then one can use high order numerical integration (e.g., tensor product Gauss-Legendre quadrature) to estimate $\mC$ in \eqref{eq:C}. In~\cite{constantine2014computing}, we propose and analyze the following Monte Carlo method for estimating $\mC$'s eigenpairs. For $i=1,\dots,M$, draw $\vx_i$ independently according to $\rho(\vx)$. Compute $\nabla f_i = \nabla f(\vx_i)$. Then estimate
\begin{equation}
\label{eq:Capprox}
\mC \;\approx\; \hmC 
\;=\; \frac{1}{M}\sum_{i=1}^M \nabla f_i\,\nabla f_i^T
\;=\; \hmW\hLambda\hmW^T.
\end{equation}
In~\cite{constantine2014computing}, we study how large $M$ must be so that (i) the estimated eigenvalues $\hLambda$ are close to the true eigenvalues $\Lambda$ and (ii) the subspace distance $\varepsilon$ defined as 
\begin{equation}
\label{eq:subdist}
\varepsilon \;=\; \|\mW_1\mW_1^T - \hmW_1\hmW_1^T\|,
\end{equation}
where $\hmW_1$ contains the first $n$ columns of $\hmW$, is well-behaved. The main result from~\cite{constantine2014computing} is Corollary 3.7, which shows that 
\begin{equation}
\label{eq:gap}
\varepsilon \;\leq\; \frac{4\,\lambda_1\,\delta}{\lambda_n - \lambda_{n+1}},
\end{equation}
where $\delta$ is a user-specified tolerance on the relative eigenvalue error that determines the minimum number of required samples. Equation \eqref{eq:gap} shows that a large eigenvalue gap $\lambda_n-\lambda_{n+1}$ implies that the active subspace can be accurately estimated with Monte Carlo. Therefore, the practical heuristic is to choose the dimension $n$ of the active subspace according to the largest gap in the eigenvalues.

\section{Sensitivity metrics from active subspaces}
\label{sec:safromas}

\noindent We propose two sensitivity metrics derived from the eigenvectors $\mW$ and eigenvalues $\Lambda$ from \eqref{eq:C}. The first metric uses the components of the first eigenvector (i.e., the eigenvector associated with the largest eigenvalue). The second metric uses a linear combination of the squared eigenvector components weighted by the eigenvalues, which we call \emph{activity scores}. 

\subsection{The first eigenvector}
\label{sec:eigenvec1}

\noindent In~\cite[Section 4.4]{constantine2015discovering}, we used the first eigenvector from the matrix $\hmC$ in \eqref{eq:Capprox} to identify the most important parameters in a single-diode model of a photovoltaic solar cell. For this particular model, the magnitudes of the eigenvector components gave the same ranking as the total sensitivity indices (see section \ref{sec:variance}). However, these two metrics measure different characteristics of the function. The total sensitivity indices measure the proportion of the variance attributed to each parameter while the eigenvector identifies an important direction in the parameter space, where ``importance'' is measured as \eqref{eq:eigenval}. In words, perturbing $\vx$ along $\vw_1$ changes $f(\vx)$ the most, on average.

Like the linear model coefficients in section \ref{sec:regression_coefficients}, the eigenvector components are not all positive, and the relative difference in signs (since eigenvectors are only unique up to a sign) provides insight into the relationship between inputs and outputs. For example, if two eigenvector components are equal in magnitude but differ in sign, then we expect $f$ to increase when one parameter is positively perturbed but decrease when the other is positively perturbed. This is similar to the gradient-based metrics described by Morris~\cite{Morris1991} and Kucherenko~\cite{kucherenko2009monte}. However, the trouble of averaging a gradient that might change signs does not effect the eigenvector. In the numerical examples in section \ref{sec:exp}, we compare the first eigenvector's components to other sensitivity metrics, but we do not provide a theoretical comparison. We use the bootstrap to estimate a standard error in the eigenvector components, though one must be careful to keep the signs consistent when computing the bootstrap replicates. Section 7.2 of~\cite{efron1994introduction} shows an example of using the bootstrap to estimate standard errors in eigenvector components. We denote the first eigenvector by $\vw_1=[w_{1,1},\dots,w_{m,1}]^T$ and the bootstrap standard error of its $i$th component by $\widehat{se}_{w_{i,1}}$.

\subsection{Activity scores}
\label{sec:activityscores}

\noindent We develop a more comprehensive and justifiable sensitivity metric from the eigenpairs in \eqref{eq:C}. We propose the following metric for global sensitivity analysis. Define $\alpha_i$ as
\begin{equation}
\label{eq:actscore}
\alpha_i \;=\; \alpha_i(n) \;=\; \sum_{j=1}^n \lambda_j\,w_{i,j}^2,
\qquad i=1,\dots,m.
\end{equation}
We call $\alpha_i(n)$ the \emph{activity score} for the $i$th parameter, and we use these numbers to rank the importance of a model's inputs. 

The active subspace's construction gives insight into the activity scores' interpretation. The eigenvector $\vw_1$ identifies the most important direction in the parameter space in the following sense: perturbing $\vx$ along $\vw_1$ changes $f$ more, on average, than perturbing $\vx$ orthogonal to $\vw_1$; see \eqref{eq:eigenval}. The components of $\vw_1$ measure the relative change in each component of $\vx$ along this most important direction, so they impart significance to each component of $\vx$. The second most important direction is the eigenvector $\vw_2$, and the relative importance of $\vw_2$ is measured by the difference between the eigenvalues $\lambda_1$ and $\lambda_2$. For example, if $\lambda_1\gg\lambda_2$, then $f(\vx)$ has a one-dimensional active subspace, and the importance of $\vx$'s components is captured in $\vw_1$'s components. Therefore, to construct the global sensitivity metric, it is reasonable to scale each eigenvector by its corresponding eigenvalue. Squaring each eigenvector component in \eqref{eq:actscore} removes information provided by the signs. But the resulting metric is much easier to compare to existing sensitivity metrics. 

\subsection{Comparison to existing metrics}
\label{sec:compare}

\noindent We compare the activity scores to the derivative-based metrics and the total sensitivity indices reviewed in Section \ref{sec:background}. 

\begin{theorem}
\label{thm:ascompare1}
The activity scores are bounded above by the derivative-based metrics,
\begin{equation}
\label{eq:as_dbsm_bound}
\alpha_i(n) \;\leq\; \nu_i, \qquad i=1,\dots,m,
\end{equation}
where $\nu_i$ is from \eqref{eq:avgsqder}. The inequality becomes an equality when $n=m$. 
\end{theorem}

\begin{pf}
Note that the derivative-based metrics are the diagonal elements of the matrix $\mC$ in \eqref{eq:C}. Let $\vnu=[\nu_1,\dots,\nu_m]^T$. Then
\begin{equation}
\vnu \;=\; \diag{\mC} \;=\; \diag{\mW\Lambda\mW^T}. 
\end{equation}
So
\begin{equation}
\nu_i \;=\; \sum_{j=1}^m \lambda_j\,w_{i,j}^2 \;=\; \alpha_i(m), \qquad i=1,\dots,m,
\end{equation}
where $w_{i,j}$ is the $(i,j)$ element of $\mW$. This proves the equality statement. To see the inequality, note that $w_{i,j}^2\geq 0$ and $\lambda_j\geq 0$, so
\begin{equation}
\alpha_{i}(n) \;=\; \sum_{j=1}^n \lambda_j\,w_{i,j}^2
\;\leq\; \sum_{j=1}^m \lambda_j\,w_{i,j}^2
\;=\; \nu_i,
\end{equation}
as required.
\end{pf}

Theorem \ref{thm:ascompare1} shows that the activity score can be interpreted as a truncation of the derivative-based metric $\nu_i$ from \eqref{eq:avgsqder}. In many practical situations, the activity scores and the derivative-based metrics give comparable rankings. However, it is possible to construct cases where the rankings differ. Consider the quadratic function $f(x_1,x_2) = (x_1^2+x_2^2)/2$. For this case, the matrix $\mC$ from \eqref{eq:C} is 
\begin{equation}
\label{eq:Cexample}
\mC \;=\; \bmat{1/3 & 0 \\ 0 & 1/3}.
\end{equation}
The derivative-based metrics are the diagonal elements of $\mC$---both $1/3$ implying equal importance in the two inputs. However, the activity scores with $n=1$ are $[1/3, 0]$, which implies that the second variable is not important at all. This case illustrates that the activity scores are appropriate when there is a gap between the eigenvalues of $\mC$ (i.e., $f$ admits an active subspace), and $n$ is chosen according to the gap. In the case of \eqref{eq:Cexample}, there is no gap between the two eigenvalues, so it is not possible to distinguish a one-dimensional eigenspace that defines a one-dimensional active subspace.

Sobol' and Kucherenko connect the total sensitivity indices to the derivative-based metrics in \cite{Sobol2009}; we use~\cite[Theorem 2]{Sobol2009} to relate the total sensitivity indices to the activity scores.

\begin{theorem}
\label{thm:ascompare2}
The total sensitivity index $\tau_i$ from \eqref{eq:saltelli_tsi} is bounded by
\begin{equation}
\tau_i \;\leq\; \frac{1}{4\pi^2\,V}\,(\alpha_i(n) + \lambda_{n+1}),
\qquad i=1,\dots,m,
\end{equation}
where $V=\Var{f}$, $\alpha_i(n)$ is from \eqref{eq:actscore}, and $\lambda_{n+1}$ is the $n+1$st eigenvalue in \eqref{eq:C}.
\end{theorem}

\begin{pf}
Theorem 2 from \cite{Sobol2009} shows 
\begin{equation}
\label{eq:ineq1}
\tau_i \;\leq \;\frac{1}{4\pi^2\,V}\,\nu_i,
\end{equation}
for $i=1,\dots,m$. Note that an additional factor of $1/4$ appears in this case due to scaling the domain and the partial derivatives to the hypercube $[-1,1]^m$ with a weight function $2^{-m}$. Using Theorem \ref{thm:ascompare1}, 
\begin{equation}
\label{eq:ineq2}
\nu_i \;=\; \alpha_i(n) + \sum_{j=n+1}^m \lambda_j w_{i,j}^2
\;\leq\; \alpha_i(n) + \lambda_{n+1} \sum_{j=n+1}^m w_{i,j}^2
\;\leq\; \alpha_i(n) + \lambda_{n+1}.
\end{equation}
The last inequality follows since the rows of $\mW$ have norm 1. Combining \eqref{eq:ineq1} with \eqref{eq:ineq2} completes the proof.
\end{pf}

Theorem \ref{thm:ascompare2} tells us that if $\lambda_{n+1}$ is small---which often indicates an $n$-dimensional active subspace---and $\alpha_i(n)$ is small, then the corresponding component of the total sensitivity index will also be small. However, Sobol and Kucherenko show an example of a function in~\cite[Section 7]{Sobol2009} where the total sensitivity indices and the derivative-based global sensitivity measures identify different sets of variables as important. The situation is similar for the activity scores. This should not be surprising since different sensitivity metrics measure different characteristics of the function. The derivative-based metrics measure the average response to small perturbations in the inputs; the activity scores are comparable when the function admits an active subspace. In contrast, the total sensitivity indices measure the variance attributable to each parameter. In many nicely behaved functions derived from practical engineering models (such as the examples in section \ref{sec:exp}), all metrics are consistent. However, it is possible to construct functions where the metrics induce different rankings. The appropriate choice of metric depends on the application. 

If the dimension of $\vx$ is sufficiently small and the derivatives of $f$ are sufficiently smooth, then one may estimate the integrals defining the matrix $\mC$ in \eqref{eq:C} using a high order numerical quadrature rule, such as tensor product Gauss-Legendre quadrature. Then the eigenvalues are computed from the numerically integrated matrix; we use this approach to compute reference values in section \ref{sec:exp}. 

When these conditions are not satisfied (e.g., if the model contains more than a handful of input parameters), which is often the case in complex simulation models, one may get coarser estimates of the activity scores $\alpha_i$ from the eigenpairs $\hLambda$, $\hmW$ of the Monte Carlo estimate $\hmC$ in \eqref{eq:Capprox}. Our previous work studied the approximation properties of the Monte Carlo-based estimates of the eigenvalues and the subspaces~\cite{constantine2014computing}. Unfortunately, those results do not directly translate to a priori error measures on the estimated activity scores. We are currently working to find a lower bound on the number $M$ of samples in \eqref{eq:Capprox} such that the activity scores are accurate to within a user-specified tolerance---similar to~\cite[Corollary 3.3]{constantine2014computing}, which bounds the number of samples needed for accurate eigenvalue estimation. We conjecture that this number may be smaller than the number of samples needed for accurate Monte Carlo estimates of the derivative-based metrics $\nu_i$ from \eqref{eq:avgsqder}, but such analysis is beyond the scope of the current paper. 

Since the Monte Carlo-based estimates of the eigenpairs cannot be interpreted as sums of independent random variables, we cannot compute central limit theorem-based standard errors. This limitation extends to the activity scores derived from the eigenpairs. Therefore, we use the bootstrap to estimate the standard error in the Monte Carlo-based estimates of the activity scores. The bootstrap algorithm is described in detail by Efron and Tibshirani~\cite[Algorithm 6.1]{efron1994introduction}. We denote the bootstrap standard error of the activity score $\alpha_i$ by $\widehat{se}_{\alpha_i}$.  In section \ref{sec:exp}, we study the bootstrap standard error's behavior as the number $M$ of samples increases in the Monte Carlo-based estimates.

\section{Numerical experiments}
\label{sec:exp}

\noindent The following numerical experiments support the claim that the proposed sensitivity metrics---namely, the first eigenvector from Section \ref{sec:eigenvec1} and the activity scores from Section \ref{sec:activityscores}---are consistent with standard sensitivity metrics. We test two models from Bingham's Virtual Library of Simulation Experiments~\cite{vlse}: (i) a nonlinear, algebraic model of a cylindrical piston with $m=7$ physical input parameters and (ii) a nonlinear, algebraic model of a transformerless push-pull circuit with $m=6$ physical input parameters. We chose these two models for two reasons. First, the algebraic expressions allow calculation of analytic gradients that are straightforward to implement, thus enabling the gradient-based metrics. Second, the models' fast numerical evaluation and low parameter space dimension allow us to compute accurate Gauss quadrature-based reference values for all quantities. Moreover, the fast evaluation enables a thorough error and convergence study for the Monte Carlo estimates of all metrics. In total, the numerical experiments used several million evaluations of each model, which would not have been possible with a more complex, expensive simulation model. We emphasize that the goal is not to study these models but to compare the sensitivity metrics. 

The parameter space in each model is a hyperrectangle defined by independent ranges on each parameter. For our studies, we normalize the parameter space to the $[-1,1]$ hypercube and scale the partial derivatives appropriately. Consequently, the sensitivity metrics measure sensitivity for the normalized model. We note that the parameter ranges are part of the definition of the model; the sensitivity metrics and active subspace computations depend on the given ranges.

For each model, we compute the active subspaces by (i) estimating the matrix $\mC$ from \eqref{eq:C} using tensor product Gauss-Legendre quadrature with 7 points per dimension and (ii) computing the eigenpairs of the quadrature-based estimate. We performed a convergence study for the quadrature-based estimates, and we found that with 7 points per dimension (i.e., $7^m$ total points), the quadrature-based estimates of the eigenvalues are accurate to 10 digits. We plot the eigenvalues on a logarithmic scale, and we find that the eigenvalues---in both models---decay spectrally, which is comparable to the real applications we have studied~\cite{Constantine2014,Lukaczyk2014,constantine2015discovering,Bauerheim2016}. We also plot the components of the first two eigenvectors, which define the first two active variables. Using the eigenvectors, we plot the one- and two-dimensional summary plots for 500 random samples of the input/output pairs. The summary plots show the relationship between the first two active variables (i.e., the two most important linear combinations of the model inputs, where the weights of the linear combinations are the first two eigenvectors of $\mC$) and the output quantity of interest. If a simple univariate relationship is present in the one-dimensional summary plot, then the output can be modeled with the low-dimensional approximation from \eqref{eq:g}. The univariate function $g$ can be fit from the samples of the active variable / output pairs; see~\cite[Chapter 4]{asm2015} for more details. If there is no apparent univariate relationship in the one-dimensional summary plot, then the two dimensional summary plot can be used similarly. If there is no apparent relationship in the two-dimensional summary plot and the eigenvalues suggest the active subspace of dimension greater than 2 is most dominant, then the summary plots are not especially helpful in revealing the low-dimensional structure in the map from inputs to outputs. Summary plots were developed in the context of sufficient dimension reduction for statistical regression. See~\cite{cook2009regression} for details in the regression context, and see~\cite[Chapter 1]{asm2015} for details in the active subspaces context. We discuss the conclusions one may draw about the active subspace from these plots before studying the sensitivity metrics. 

For each of the two models, we examine the following five sensitivity metrics:
\begin{enumerate}
\item the Sobol' total sensitivity index $\tau_i$ from Section \ref{sec:variance},
\item the derivative-based global sensitivity measures $\nu_i$ from Section \ref{sec:derivative},
\item the linear model coefficients $\beta_i$ from Section \ref{sec:regression_coefficients},
\item the first eigenvector components $w_{i,1}$ from Section \ref{sec:eigenvec1},
\item and activity scores $\alpha_i(n)$, where $n$ is the dimension of the active subspace, from Section \ref{sec:activityscores}.
\end{enumerate}
The subscript $i$ gives the index of the corresponding parameter. For each sensitivity metric, we compute
\begin{itemize}
\item a tensor product Gauss-Legendre quadrature-based reference value with 7 points per dimension (i.e., $7^m$ total points)\footnote{With 7 points per dimension, the reference values contain 10 digits of accuracy in all cases.},
\item a Monte Carlo estimate with $M$ samples, where $M$ ranges between 50 and 50000,
\item a relative error between the Monte Carlo estimate and the quadrature-based reference value, 
\item and the (bootstrap) standard error for each Monte Carlo estimate.
\end{itemize}
The number $M$ of evaluations of $f(\vx)$ is a measure of the sensitivity metrics' cost. For the total sensitivity indices, we choose $M'$ such that $M=(m+1)M'$, and we use $M'$ as the number of samples in $\vf_{\mA}$ and each $\vf_{\mB_i}$ from Section \ref{sec:variance}. This keeps the comparison fair, but the number of samples used to compute the integrals is effectively less than for the other metrics. Consequently, we expect the error and standard error to be larger for the total sensitivity indices. For the derivative-based metrics, we count each gradient sample as one evaluation, which is the case for simulation models with algorithmic differentiation capabilities.

We scale the error and standard error in the following manner. Let $\gamma_i$ be the reference value for one of the five sensitivity metrics listed above, and let $\hat{\gamma}_i$ be its Monte Carlo estimate. We compute the relative error in the Monte Carlo estimate as
\begin{equation}
\label{eq:relerr}
\frac{|\gamma_i - \hat{\gamma}_i|}{
\max\limits_{1\leq j\leq m}\;|\gamma_j|
}.
\end{equation}
Scaling error by the magnitude of the largest-in-magnitude metric avoids issues of dividing by very small sensitivity metrics. Since these errors are random due to the random sampling, we average the errors, as a function of the number $M$ of samples, over 10 independent trials. We divide the (bootstrap) standard errors by the magnitude of the largest-in-magnitude Monte Carlo estimate of the sensitivity metric; in other words, we do not use the reference values when computing the standard errors as would be done in practice without reference values. We average the standard error estimates over 10 independent trials to reduce the effects of randomness. For the active subspace-based activity scores, we should how the values change as a function of the dimension $n$ of the active subspace. Recall that Theorem \ref{thm:ascompare1} shows that the activity scores become the derivative-based global sensitivity metrics when the dimension $n$ of the active subspace reaches the dimension $m$ of the parameter space. 

Finally, we show that the rankings on the input variables' importance induced by the five sensitivity metrics are consistent for both models. We emphasize that this is not case in general. Since the various metrics are derived from different features of output (e.g., averaged derivatives versus variances), it is possible to devise functions where the metrics are not consistent; see, e.g., \cite[Section 7]{Sobol2009}. The purpose of comparing rankings is to show that the metrics derived from the active subspace are consistent with more standard global sensitivity metrics in an engineering inspired model, which supports the analysis from Section \ref{sec:compare}. 

The Matlab scripts for all computations and plots are available at \url{https://bitbucket.org/paulcon/global-sensitivity-metrics-from-active-subspaces}. All computations were run in Matlab 2015b on a MacBook Pro with a dual core processor and 16GB of memory. 

\subsection{Piston model}
\label{sec:piston}

\noindent The first model is an algebraic expression for a cylindrical piston. Details of the model and associated Matlab code are available at \url{http://www.sfu.ca/~ssurjano/piston.html}. The piston model appears in~\cite{ben2007modeling,moon2010design} as a test model for statistical screening. The quantity of interest, $t$, is the time in seconds it takes the piston to complete one cycle. This time depends on $m=7$ physical parameters according to the following nonlinear expressions:
\begin{align}
\label{eq:piston}
t &= 2\pi 
\sqrt{\frac{M}{k+S^2\frac{P_0 V_0}{T_0}\frac{T_a}{V^2}}},\\
V &= \frac{S}{2k}\left(
\sqrt{A^2+4k\frac{P_0 V_0}{T_0}T_a}-A 
\right),\\
A &= P_0 S+19.62M - \frac{k V_0}{S}.
\end{align}
The input parameters' descriptions, ranges, and units are in Table \ref{tab:piston_input}. 

\begin{table}[!h]
\centering
\caption{Input parameters' descriptions, ranges, and units for the piston model \eqref{eq:piston}.}
\begin{tabular}{ l l l l l }
Parameter & Notation & Min & Max & Units\\
\hline
piston weight & $M$ & 30 & 60 & kg\\
piston surface area & $S$ & 0.005 & 0.020 & m$^2$ \\
initial gas volume & $V_0$ & 0.002 & 0.010 & m$^3$\\
spring coefficient & $k$ & 1000 & 5000 & N/m\\
atmospheric pressure & $P_0$ & 90000 & 110000 & N/m$^2$\\
ambient temperature & $T_a$ & 290 & 296 & K \\
filling gas temperature & $T_0$ & 340 & 360 & K
\end{tabular}
\label{tab:piston_input}
\end{table}

\begin{figure}[!h]
\centering
\subfloat[Eigenvalues of $\mC$]{%
\includegraphics[width=0.45\textwidth]{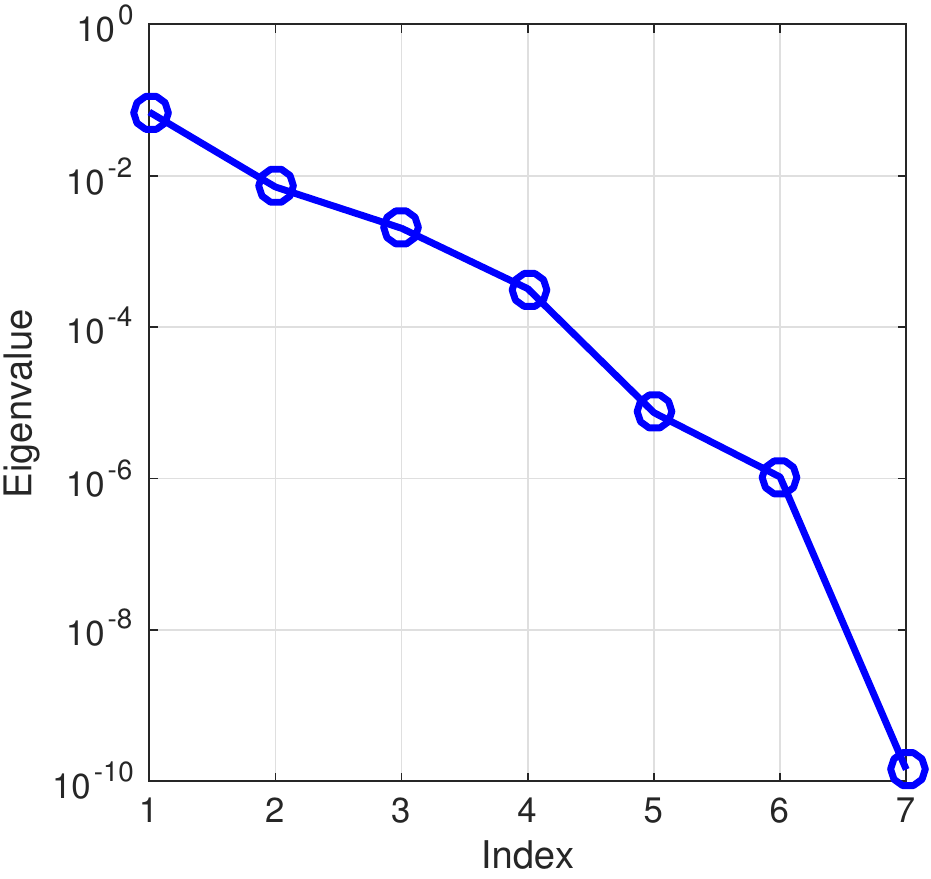}
\label{subfig:as_evals_piston}
}
\subfloat[Two eigenvectors of $\mC$]{%
\includegraphics[width=0.45\textwidth]{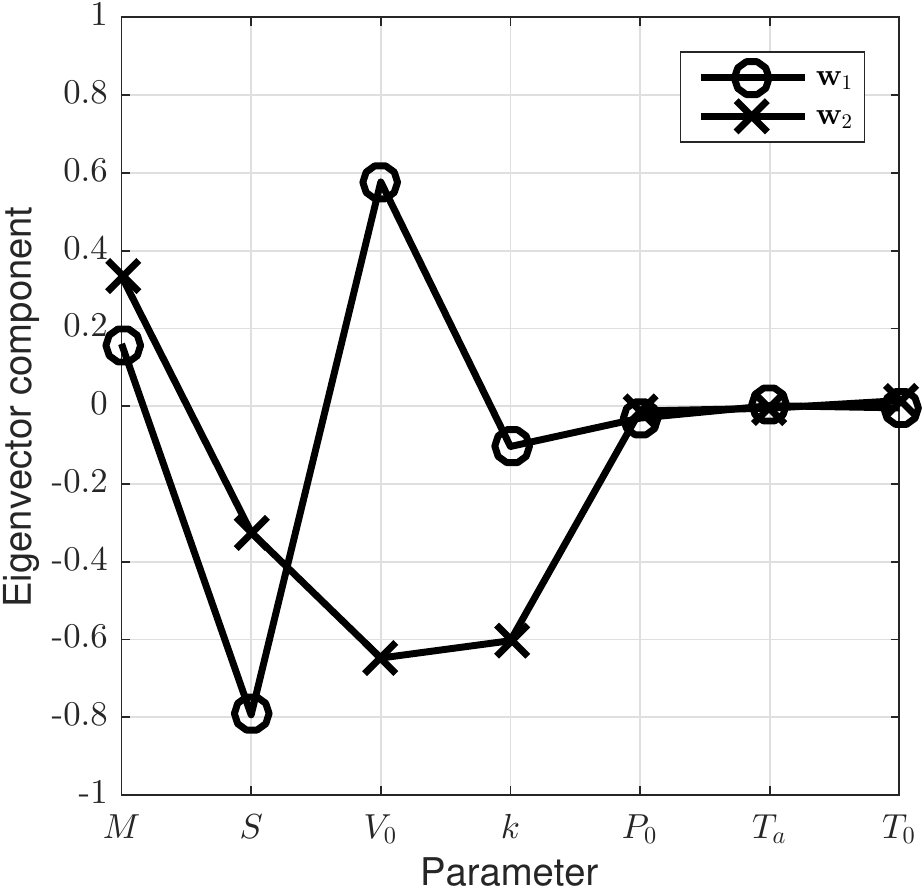}
\label{subfig:as_evecs_piston}
}\\
\subfloat[1D summary plot]{%
\includegraphics[width=0.45\textwidth]{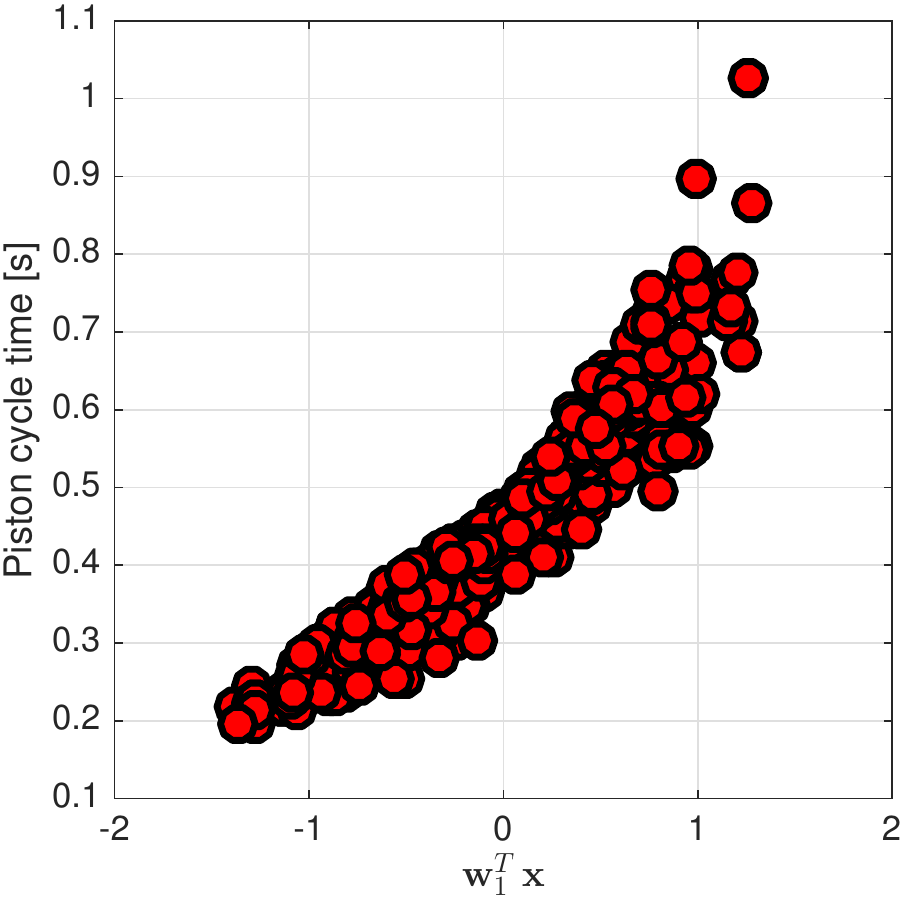}
\label{subfig:as_ss1_piston}
}
\subfloat[2D summary plot]{%
\includegraphics[width=0.45\textwidth]{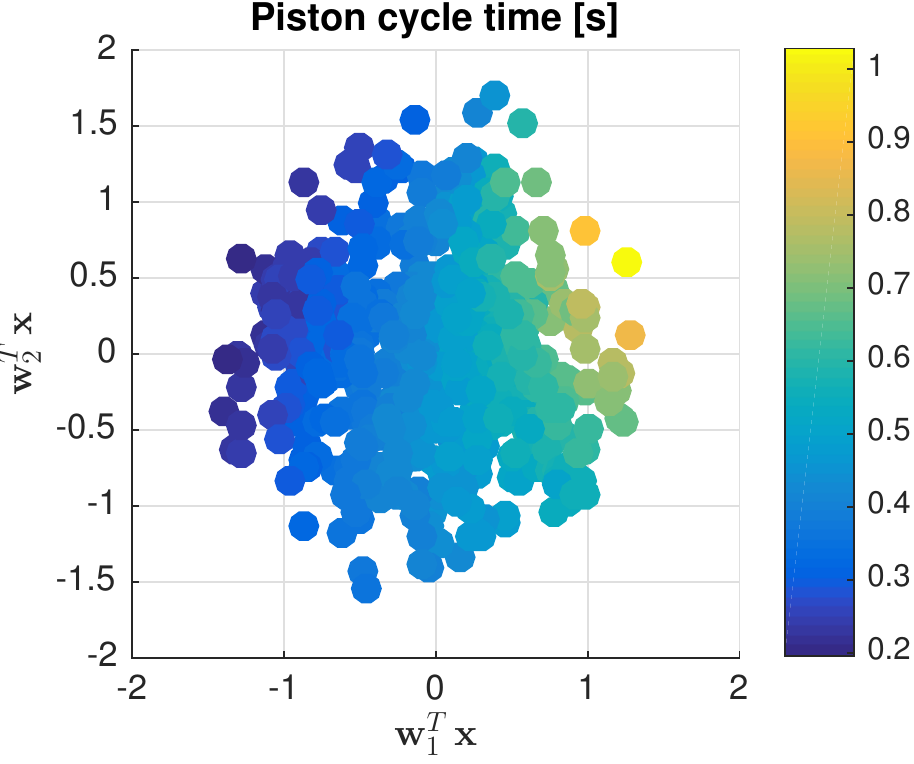}
\label{subfig:as_ss2_piston}
}
\caption{The eigenvalues in the top left subfigure provide evidence of active subspaces of dimension 1 to 6 in the piston cycle time as a function of its 7 physical input parameters from Table \ref{tab:piston_input}. Components of the first two eigenvectors are in the top right. The one- and two-dimensional summary plots in the second row elucidate the relationship between the first two active variables and the piston cycle time.}
\label{fig:as_piston}
\end{figure}

Figure \ref{subfig:as_evals_piston} shows the 7 eigenvalues, on a logarithmic scale, from the quadrature-based estimate of $\mC$ from \eqref{eq:C}. The order-of-magnitude gaps between the eigenvalues suggest that an active subspace exists for $n$ from 1 to 6. Figure \ref{subfig:as_evecs_piston} shows the components of the first two eigenvectors of $\mC$, which are used to produce the one- and two-dimensional summary plots in Figures \ref{subfig:as_ss1_piston} and \ref{subfig:as_ss2_piston}, respectively. The one-dimensional summary plot shows a gross monotonic trend in the output as a function of the first active variable. This trend might be useful in certain applications, such as determining a good starting value for a numerical optimization routine. See~\cite{constantine2014exploiting} for an example of exploiting such a relationship for estimating a range of outputs in a multiphysics scramjet model. However, there is also significant spread around a perceived univariate relationship, so an approximation of the form \eqref{eq:g} with one eigenvector may not be useful for all response surface applications. This spread is mostly likely due to changes in the output as the input is perturbed orthogonally to the first eigenvector of $\mC$. The two-dimensional summary plot captures more of the input/output relationship, where some curvature is apparent in the two-dimensional level sets. Thus, an approximation of the form \eqref{eq:g} with two eigenvectors may be more useful than with only one eigenvector. See~\cite{Lukaczyk2014} for an example of building a response surface that exploits the two-dimensional active subspace for aerospace shape design. This sort of qualitative reasoning is typical when the engineer is determining how to exploit the active subspace for a particular application.

Table \ref{tab:piston_ref} shows the reference values for each sensitivity metric to four digits. Since the eigenvalue gap is largest between the first and second eigenvalues, we choose $n=1$ when computing the activity scores to compare to the other sensitivity metrics. These are the values used to compute the error of the Monte Carlo estimates.

\begin{table}[!h]
\centering 
\caption{Reference values for piston model's sensitivity metrics computed with tensor product Gauss-Legendre quadrature. The metrics compared are the Sobol' total sensitivity index ($\tau_i$, \eqref{eq:tsi}), the derivative-based global sensitivity measure ($\nu_i$, \eqref{eq:avgsqder}), the linear model coefficients ($\beta_i$, \eqref{eq:srq}), the first eigenvector components ($w_{i,1}$, Section \ref{sec:eigenvec1}), and the activity scores ($\alpha_i$, \eqref{eq:actscore}).}
\begin{tabular}{llllll}
Parameter & $\tau_i$ &  $\nu_i$ & $\beta_i$ & $w_{i,1}$ & $\alpha_i(1)$ \\  
\hline
$M$ & 0.0509 & 0.0032 & 0.1963 & 0.1604 & 0.0018\\
$S$ & 0.5994 & 0.0449 & -0.7345 & -0.7936 & 0.0437\\
$V_0$ & 0.3528 & 0.0265 & 0.5567 & 0.5768 & 0.0231\\
$k$ & 0.0669 & 0.0040 & -0.1424 & -0.1035 & 0.0007\\
$P_0$ & 0.0013 & 0.0001 & -0.0352 & -0.0305 & 0.0001\\
$T_a$ & 0.0000 & 0.0000 & 0.0018 & 0.0015 & 0.0000\\
$T_0$ & 0.0001 & 0.0000 & -0.0051 & -0.0042 & 0.0000
\end{tabular}
\label{tab:piston_ref}
\end{table}

\begin{figure}[!h]
\centering
\subfloat[Total sensitivity index, $\tau_i$]{%
\includegraphics[width=0.32\textwidth]{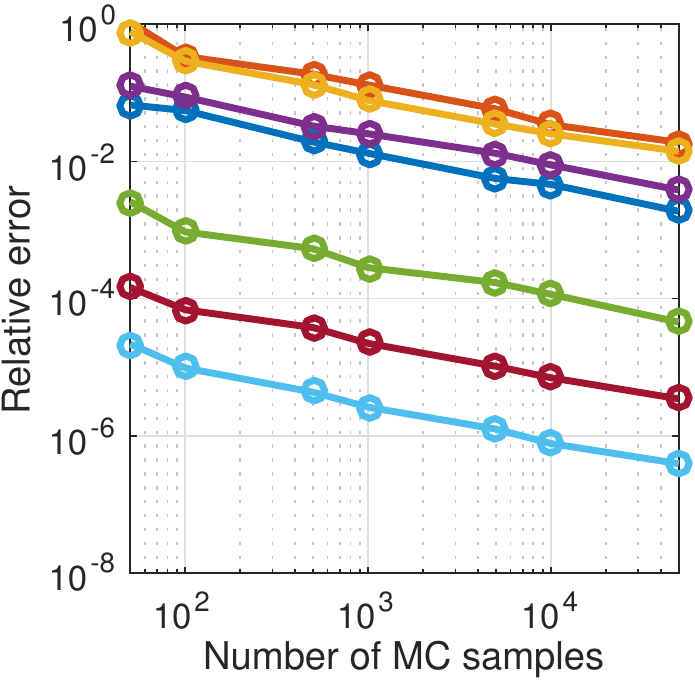}
\label{fig:err_tsi_piston}
}
\subfloat[Derivative-based metric, $\nu_i$]{%
\includegraphics[width=0.32\textwidth]{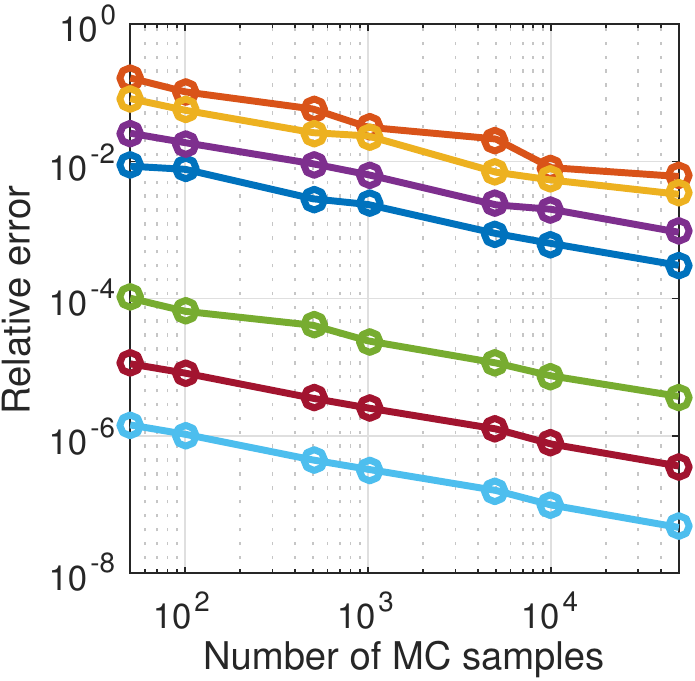}
\label{fig:err_dbsm_piston}
}
\subfloat[Linear model coefficient, $\beta_i$]{%
\includegraphics[width=0.32\textwidth]{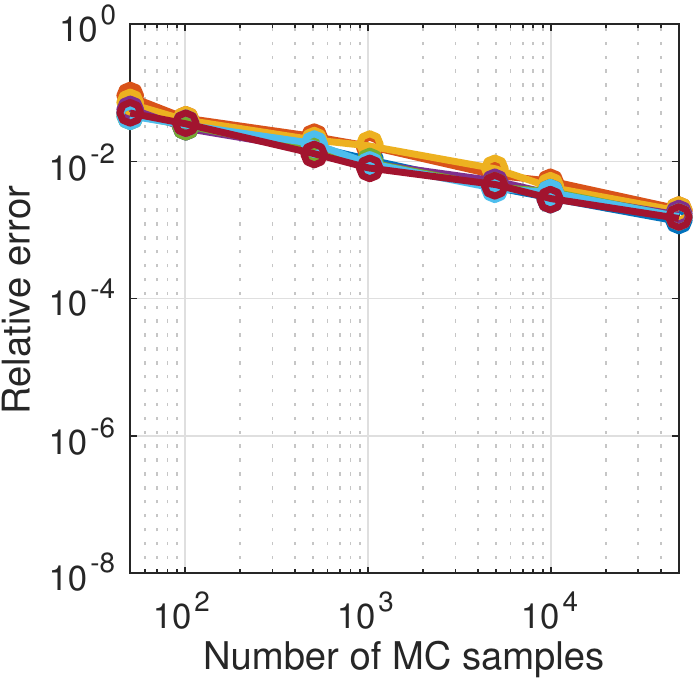}
\label{fig:err_reg_piston}
}\\
\subfloat[First eigenvector, $w_{i,1}$]{%
\includegraphics[width=0.32\textwidth]{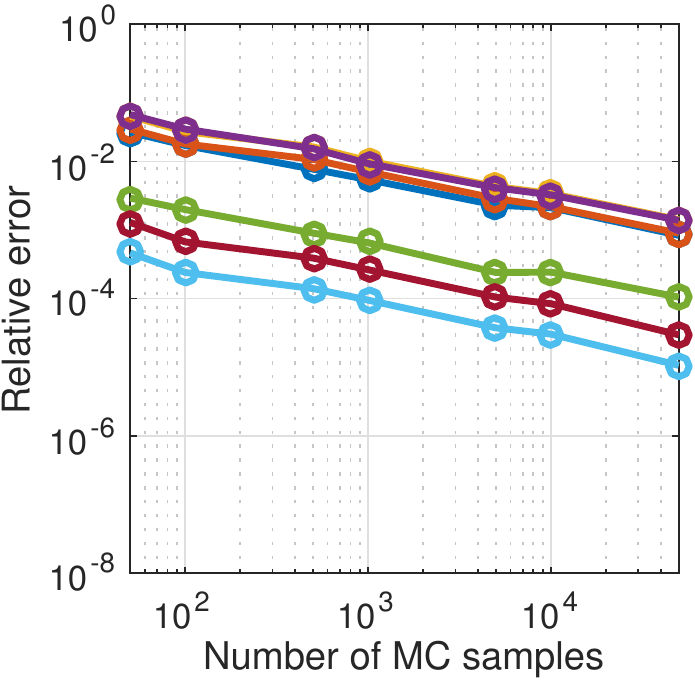}
\label{fig:err_w1_piston}
}
\subfloat[Activity score, $\alpha_i(1)$]{%
\includegraphics[width=0.32\textwidth]{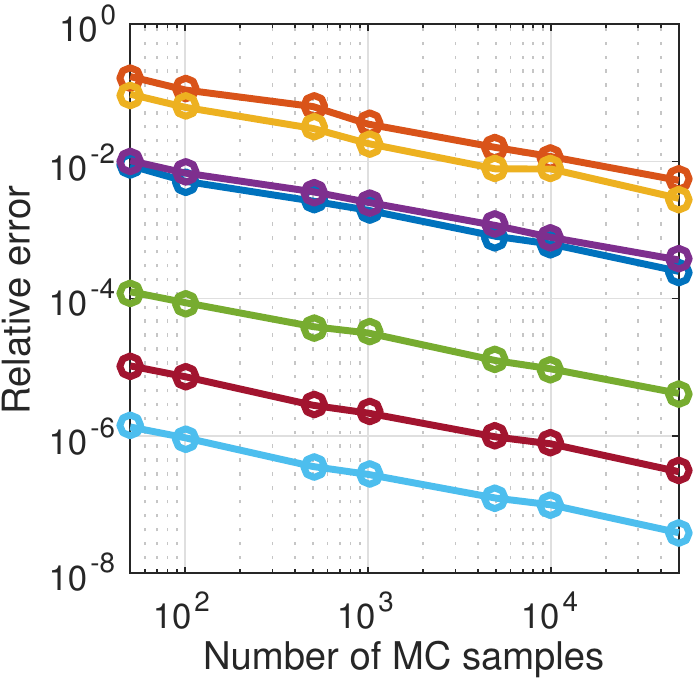}
\label{fig:err_as_piston}
}\hfil
\subfloat[Legend]{%
\includegraphics[width=0.25\textwidth]{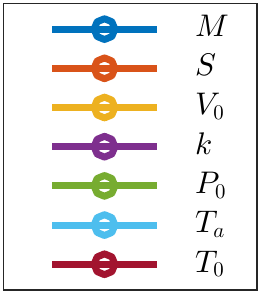}
}\hfil
\caption{
The relative error in the Monte Carlo estimates of the sensitivity metrics for the piston model \eqref{eq:piston} as a function of the number $M$ of samples, computed using the reference values in Table \ref{tab:piston_ref} and the formula in \eqref{eq:relerr}.
}
\label{fig:piston_err}
\end{figure}

Figure \ref{fig:piston_err} shows relative error in the Monte Carlo estimates of the sensitivity metrics, averaged over 10 independent trials. All errors decrease at the expected $M^{-1/2}$ convergence rate. The errors for the total sensitivity indices are larger due to how we allocate the $M$ model evaluations as described above. Across all metrics, the errors in the large metrics are relatively larger than the errors in the small metrics---with the notable exception of the linear model coefficients, which have roughly the same error. This is likely due to our choice for normalization in \eqref{eq:relerr}. 

\begin{figure}[!h]
\centering
\subfloat[Total sensitivity index, $\tau_i$]{%
\includegraphics[width=0.32\textwidth]{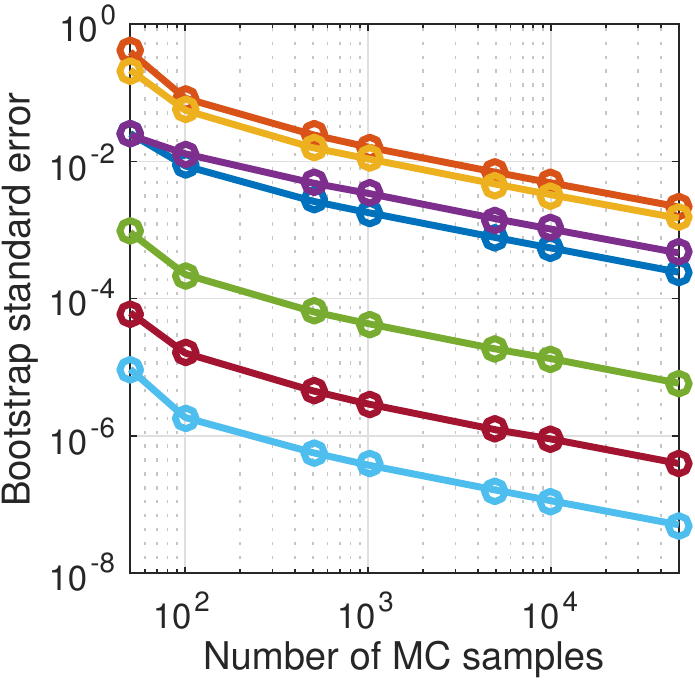}
\label{fig:stderr_tsi_piston}
}
\subfloat[Derivative-based metric, $\nu_i$]{%
\includegraphics[width=0.32\textwidth]{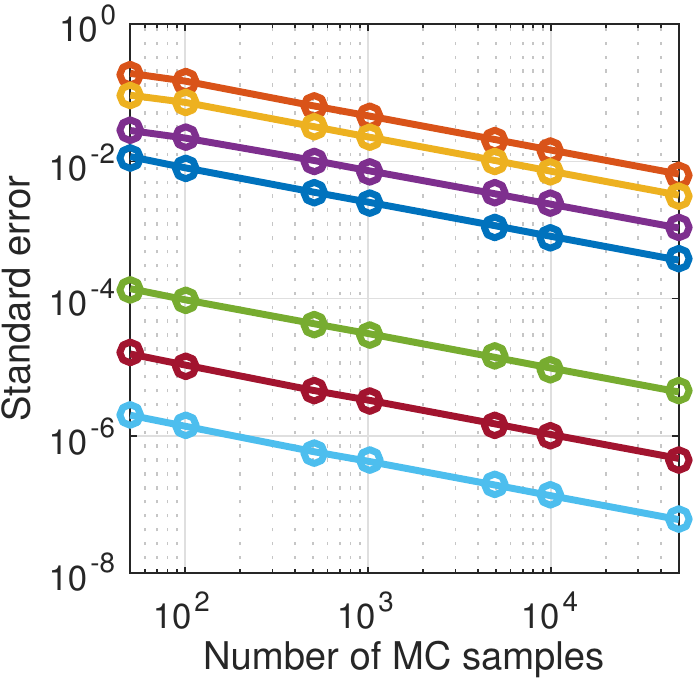}
\label{fig:stderr_dbsm_piston}
}
\subfloat[Linear model coefficient, $\beta_i$]{%
\includegraphics[width=0.32\textwidth]{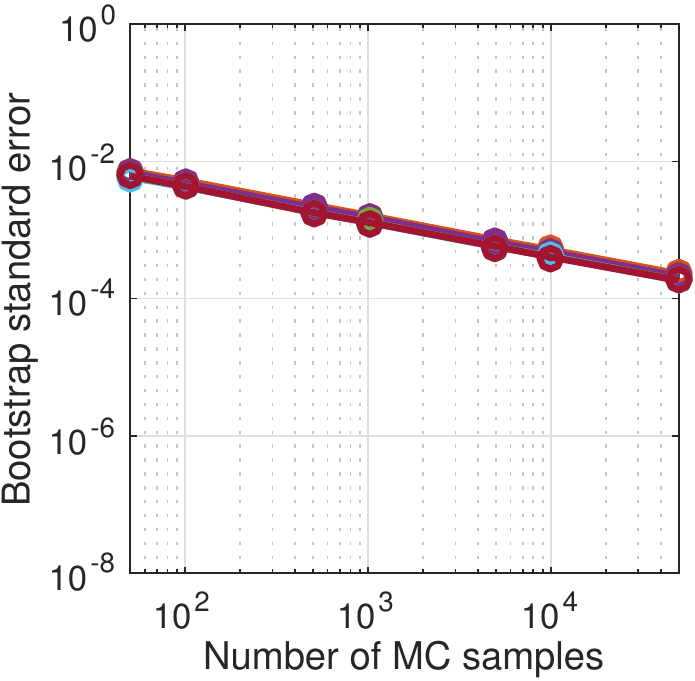}
\label{fig:stderr_reg_piston}
}\\
\subfloat[First eigenvector, $w_{i,1}$]{%
\includegraphics[width=0.32\textwidth]{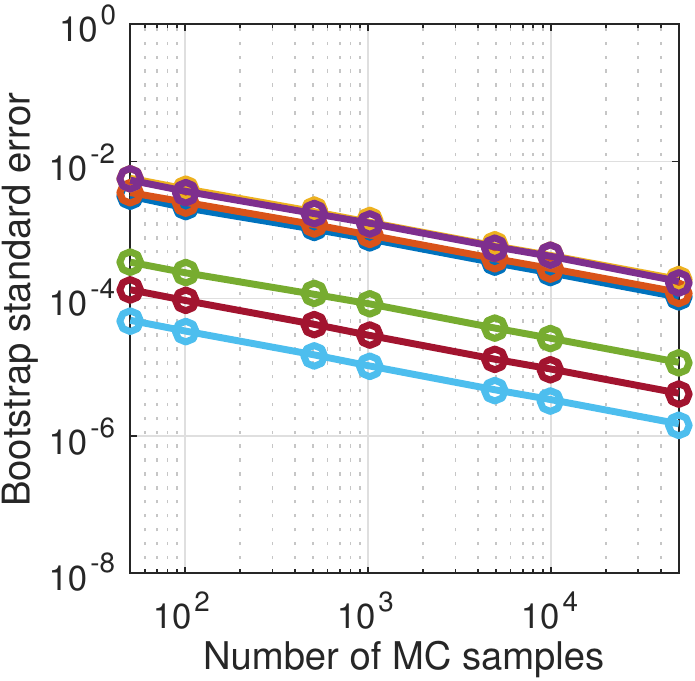}
\label{fig:stderr_w1_piston}
}
\subfloat[Activity score, $\alpha_i(1)$]{%
\includegraphics[width=0.32\textwidth]{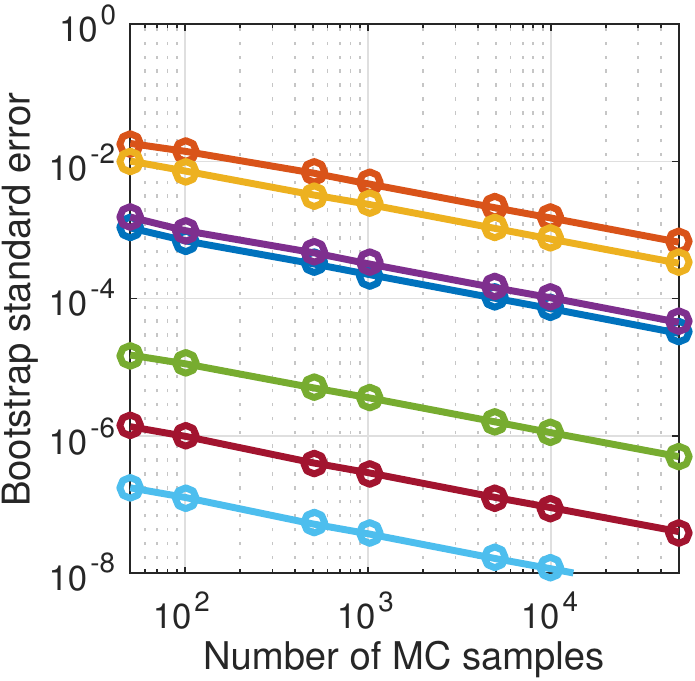}
\label{fig:stderr_as_piston}
}\hfil
\subfloat[Legend]{%
\includegraphics[width=0.25\textwidth]{figs/piston/legend.pdf}
}\hfil
\caption{
The (bootstrap) standard errors in the Monte Carlo estimates of the sensitivity metrics for the piston model \eqref{eq:piston} as a function of the number $M$ of samples.
}
\label{fig:piston_stderr}
\end{figure}

Figure \ref{fig:piston_stderr} shows the (bootstrap) standard errors for the Monte Carlo estimates. The bootstrap standard errors are generally lower than the relative errors from Figure \ref{fig:piston_err}. However, all standard errors decay at the expected rate of $M^{-1/2}$. 

\subsection{Circuit model}
\label{sec:otl}

\noindent The second model is an algebraic expression for a transformerless push-pull circuit. Details of the model and associated Matlab code are available at \url{http://www.sfu.ca/~ssurjano/otlcircuit.html}. The circuit model also appears in~\cite{ben2007modeling,moon2010design} as a test model for statistical screening. The quantity of interest, $V$, is the midpoint voltage of the circuit. This voltage depends on $m=6$ physical parameters according to the following nonlinear expressions:
\begin{align}
\label{eq:otl}
V &= 
\frac{
(12R_{b2}/(R_{b1}+R_{b2})+0.74)\beta(R_{c2}+9)
}{
\beta (R_{c2}+9)+R_f
}\\
&\quad +\;\frac{11.35 R_f}{\beta (R_{c2}+9) +R_f}
 + \frac{0.74 R_f\beta(R_{c2}+9)}{(\beta(R_{c2}+9)+R_f)R_{c1}}.
\end{align}
The input parameters' descriptions, ranges, and units are in Table \ref{tab:otl_input}. 

\begin{table}[!h]
\centering 
\caption{Input parameters' descriptions, ranges, and units for the circuit model \eqref{eq:otl}.}
\begin{tabular}{lllll}
Parameter & Notation & Min & Max & Units\\
\hline
resistance b1 & $R_{b1}$ & 50 & 150 & K-Ohms\\
resistance b2 & $R_{b2}$ & 25 & 70 & K-Ohms\\
resistance f & $R_f$ & 0.5 & 3.0 & K-Ohms\\
resistance c1 & $R_{c1}$ & 1.2 & 2.5 & K-Ohms\\
resistance c2 & $R_{c2}$ & 0.25 & 1.20 & K-Ohms\\
current gains & $\beta$ & 50 & 300 & Amperes
\end{tabular}
\label{tab:otl_input}
\end{table}

\begin{figure}[!h]
\centering
\subfloat[Eigenvalues of $\mC$]{%
\includegraphics[width=0.45\textwidth]{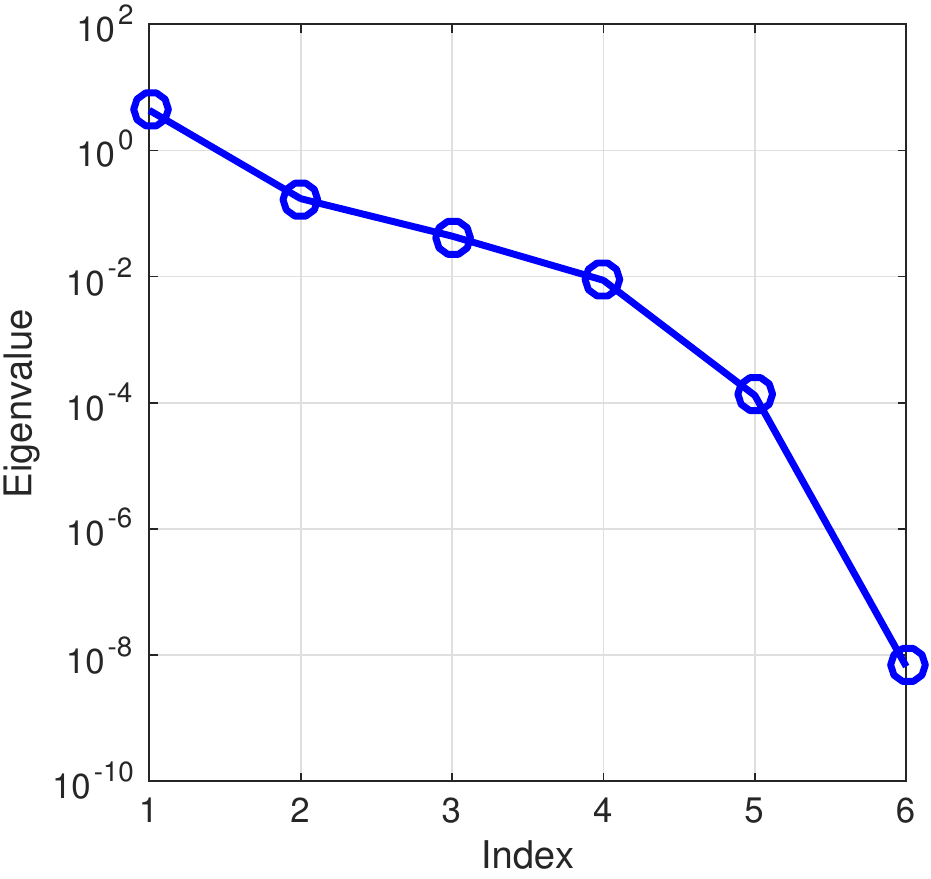}
\label{subfig:as_evals_otlcircuit}
}
\subfloat[Two eigenvectors of $\mC$]{%
\includegraphics[width=0.45\textwidth]{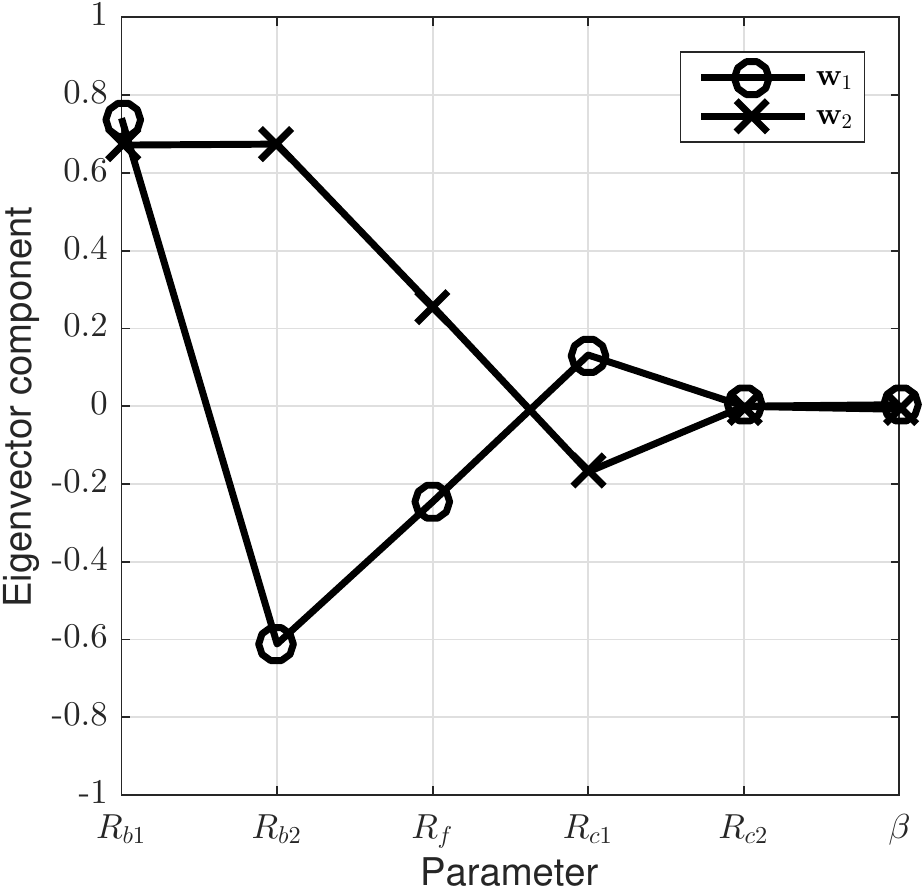}
\label{subfig:as_evecs_otlcircuit}
}\\
\subfloat[1D summary plot]{%
\includegraphics[width=0.45\textwidth]{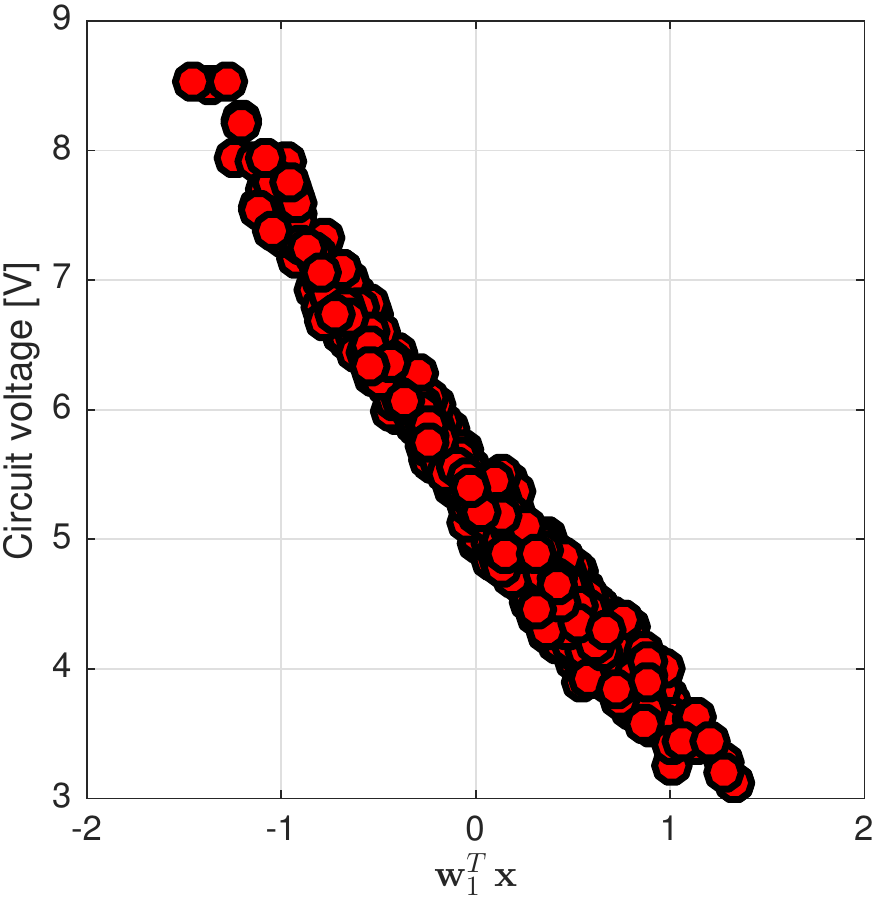}
\label{subfig:as_ss1_otlcircuit}
}
\subfloat[2D summary plot]{%
\includegraphics[width=0.45\textwidth]{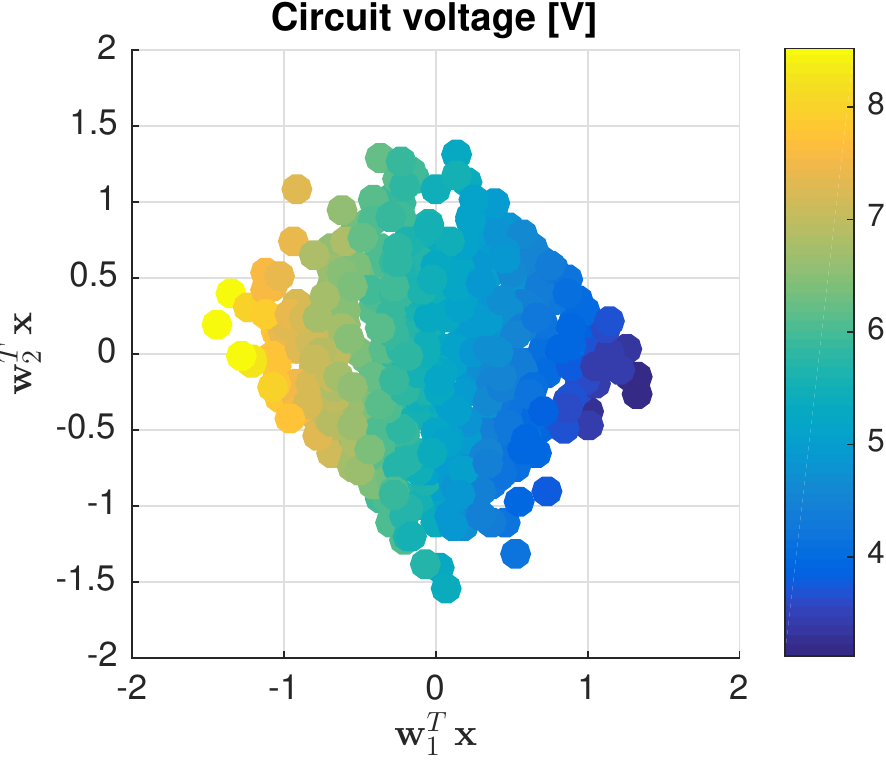}
\label{subfig:as_ss2_otlcircuit}
}
\caption{The eigenvalues in the top left subfigure provide evidence of active subspaces of dimension 1 to 5 in the circuit voltage as a function of its 6 physical input parameters from Table \ref{tab:otl_input}. Components of the first two eigenvectors are in the top right. The one- and two-dimensional summary plots in the second row elucidate the relationship between the first two active variables and the circuit voltage.}
\label{fig:as_otlcircuit}
\end{figure}

Figure \ref{subfig:as_evals_otlcircuit} shows the 6 eigenvalues, on a logarithmic scale, from the quadrature-based estimate of $\mC$ from \eqref{eq:C}. The order-of-magnitude gaps between the eigenvalues suggest that an active subspace exists for $n$ from 1 to 5. Figure \ref{subfig:as_evecs_otlcircuit} shows the components of the first two eigenvectors of $\mC$, which are used to produce the one- and two-dimensional summary plots in Figures \ref{subfig:as_ss1_otlcircuit} and \ref{subfig:as_ss2_otlcircuit}, respectively. The one-dimensional summary plot shows a nearly linear relationship in the output as a function of the first active variable. Compared to the piston model (see Figure \ref{subfig:as_ss1_piston}), the univariate relationship is tighter, and the spread away from the univariate relationship is relatively small. This suggests that, relative to the structure in the piston model, the relationship between the circuit parameters and the output voltage can be better modeled by a function---even a linear function---of the first active variable. In this case $g$ from \eqref{eq:g} may be a univariate linear function fit with least-squares; see~\cite[Chapter 4]{asm2015} for more details. This is confirmed by the two-dimensional summary plot, which has nearly linear contours from the 500 random samples. Again, this is the sort of qualitative reasoning an engineer may use to determine how to exploit the low-dimensional structure for a particular application.

Table \ref{tab:otl_ref} shows the reference values for each sensitivity metric to four digits. Similar to the piston model, since the eigenvalue gap is largest between the first and second eigenvalues, we choose $n=1$ when computing the activity scores to compare to the other sensitivity metrics. These are the values used to compute the error of the Monte Carlo estimates.

\begin{table}[!h]
\centering 
\caption{%
Reference values for circuit model's sensitivity metrics. The notation is similar to Table \ref{tab:piston_ref}. 
}
\begin{tabular}{llllll}
Parameter & $\tau_i$ & $\nu_i$ & $\beta_i$ & $w_{i,1}$ & $\alpha_i(1)$\\
\hline
$R_{b1}$ & 0.5001 & 2.4555 & -0.6925 & -0.7407 & 2.3778 \\
$R_{b2}$ & 0.4117 & 1.7039 & 0.6358 & 0.6112 & 1.6190 \\
$R_f$ & 0.0740 & 0.2902 & 0.2662 & 0.2458 & 0.2617 \\
$R_{c1}$ & 0.0218 & 0.1090 & -0.1341 & -0.1318 & 0.0752 \\
$R_{c2}$ & 0.0000 & 0.0000 & -0.0002 & -0.0002 & 0.0000 \\
$\beta$ & 0.0000 & 0.0002 & -0.0032 & -0.0039 & 0.0001
\end{tabular}
\label{tab:otl_ref}
\end{table}

\begin{figure}[!h]
\centering
\subfloat[Total sensitivity index]{%
\includegraphics[width=0.32\textwidth]{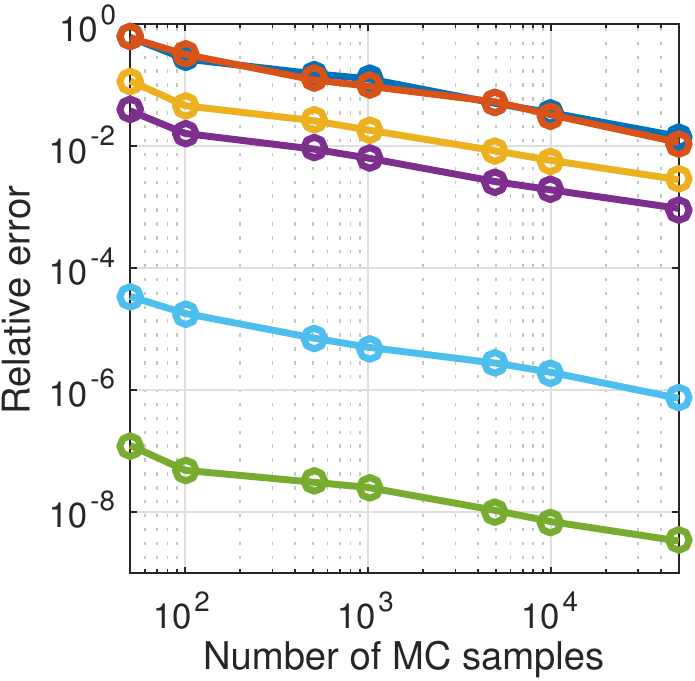}
\label{fig:err_tsi_otlcircuit}
}
\subfloat[Derivative-based metric]{%
\includegraphics[width=0.32\textwidth]{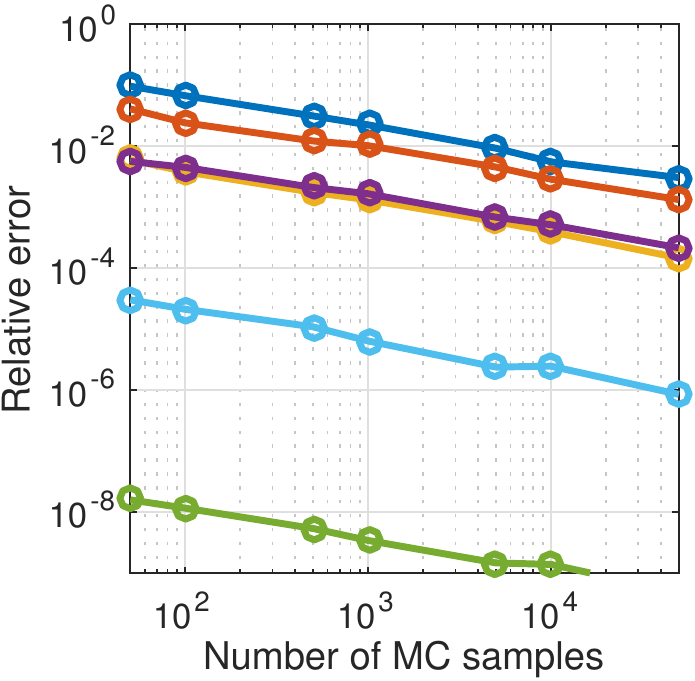}
\label{fig:err_dbsm_otlcircuit}
}
\subfloat[Linear model coefficient]{%
\includegraphics[width=0.32\textwidth]{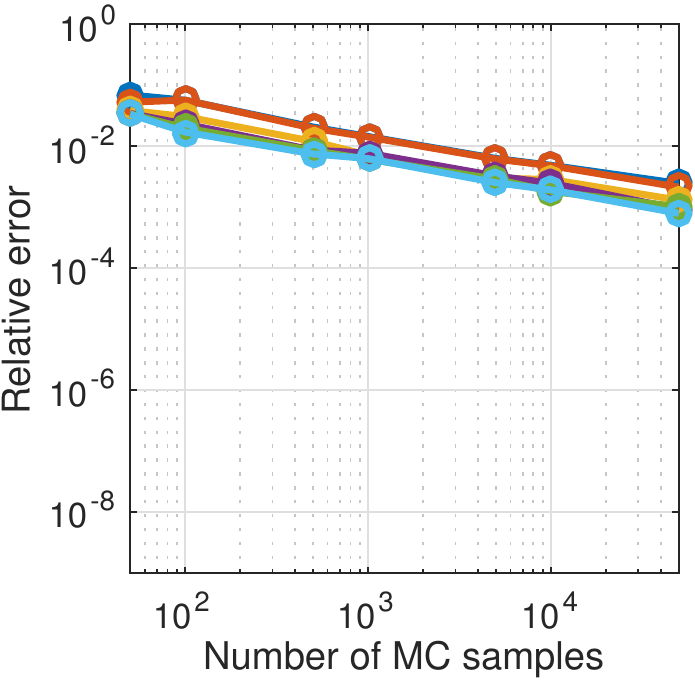}
\label{fig:err_reg_otlcircuit}
}\\
\subfloat[First eigenvector]{%
\includegraphics[width=0.32\textwidth]{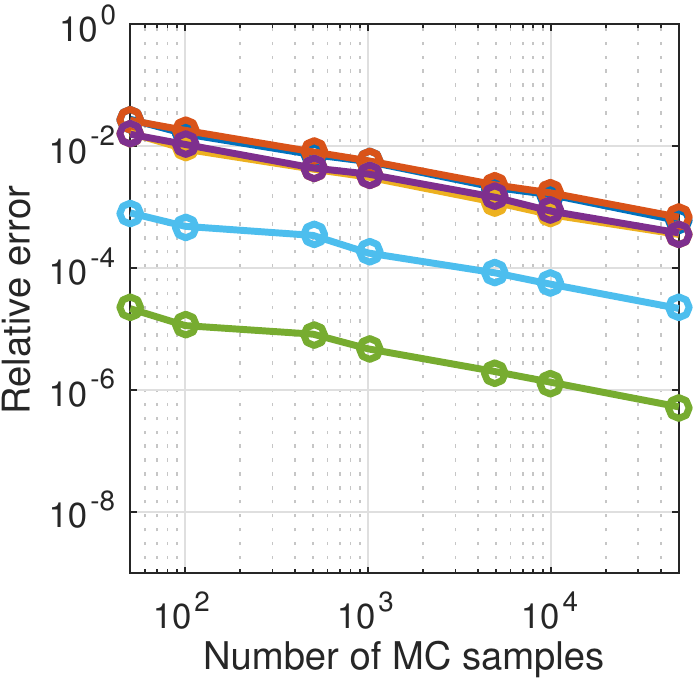}
\label{fig:err_w1_otlcircuit}
}
\subfloat[Activity score]{%
\includegraphics[width=0.32\textwidth]{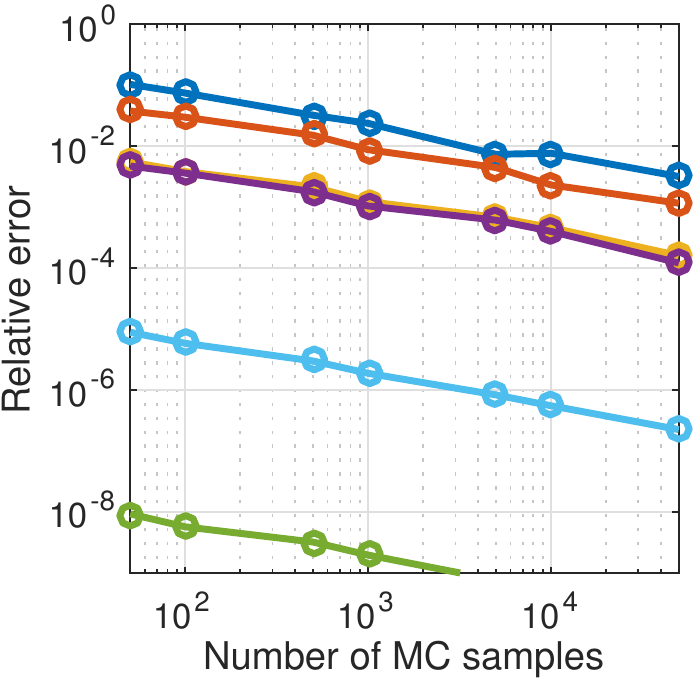}
\label{fig:err_as_otlcircuit}
}\hfil
\subfloat[Legend]{%
\includegraphics[width=0.25\textwidth]{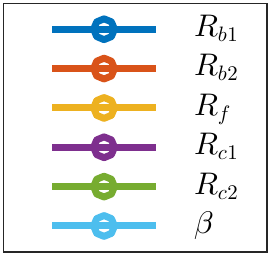}
}\hfil
\caption{%
The relative error in the Monte Carlo estimates of the sensitivity metrics for the circuit model \eqref{eq:otl} as a function of the number $M$ of samples, averaged over 10 independent trials and computed using the error definition \eqref{eq:relerr} with reference values from Table \ref{tab:otl_input}.
}
\label{fig:otl_err}
\end{figure}

Figure \ref{fig:otl_err} shows relative errors in the Monte Carlo estimates of the sensitivity metrics, averaged over 10 independent trials. All errors decrease at the expected $M^{-1/2}$ convergence rate. The errors in the total sensitivity indices are relatively large due to how we allocate the samples. Similar to the piston model, the errors in the large metrics are relatively larger than the errors in the small metrics---except for the linear model coefficients, which have roughly the same error. 

\begin{figure}[!h]
\centering
\subfloat[Total sensitivity index]{%
\includegraphics[width=0.32\textwidth]{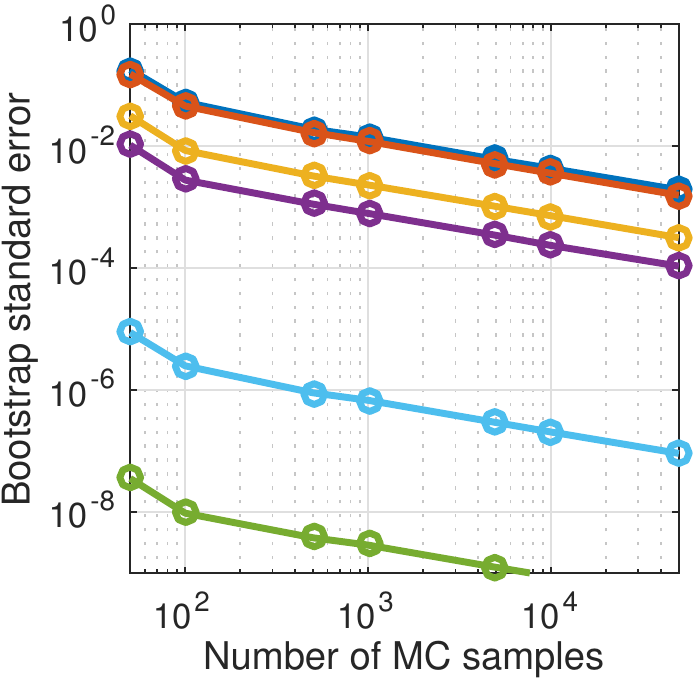}
\label{fig:stderr_tsi_otlcircuit}
}
\subfloat[Derivative-based metric]{%
\includegraphics[width=0.32\textwidth]{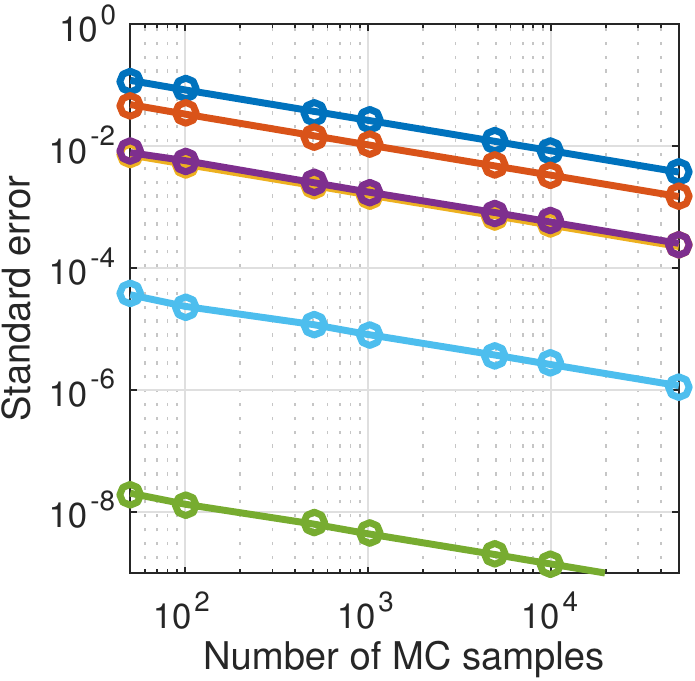}
\label{fig:stderr_dbsm_otlcircuit}
}
\subfloat[Linear model coefficient]{%
\includegraphics[width=0.32\textwidth]{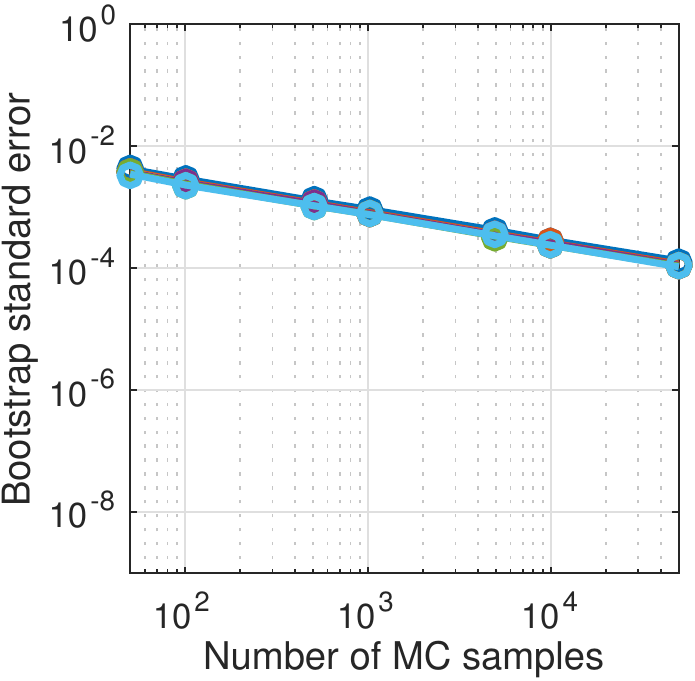}
\label{fig:stderr_reg_otlcircuit}
}\\
\subfloat[First eigenvector]{%
\includegraphics[width=0.32\textwidth]{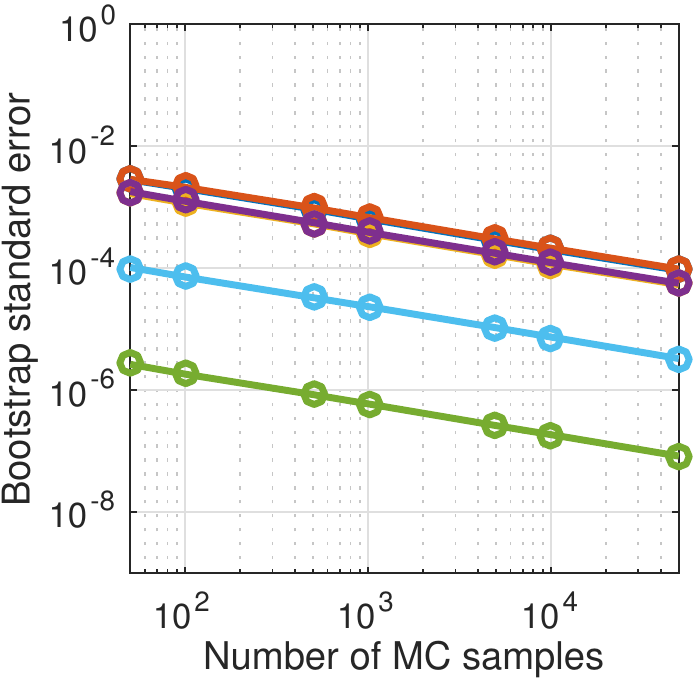}
\label{fig:stderr_w1_otlcircuit}
}
\subfloat[Activity score]{%
\includegraphics[width=0.32\textwidth]{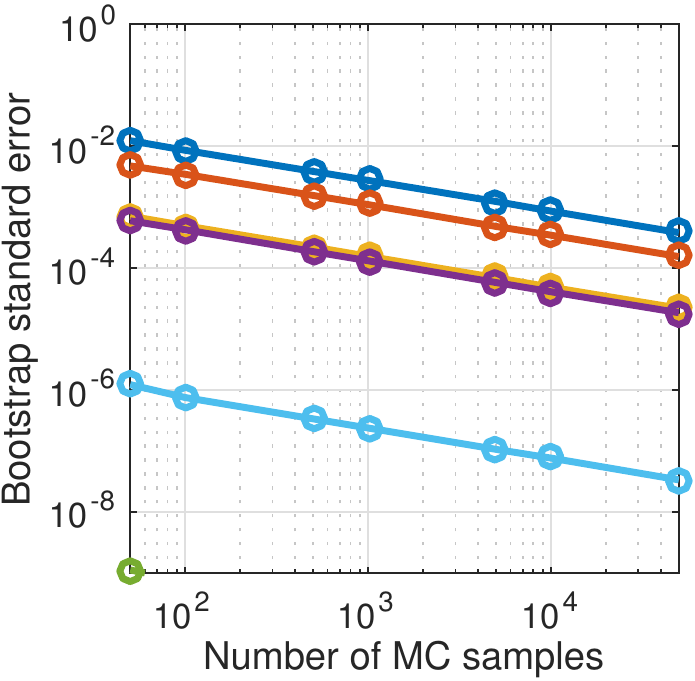}
\label{fig:stderr_as_otlcircuit}
}\hfil
\subfloat[Legend]{%
\includegraphics[width=0.25\textwidth]{figs/otlcircuit/legend.pdf}
}\hfil
\caption{The (bootstrap) standard errors for the Monte Carlo estimates as a function of the number $M$ of samples averaged over 10 independent trials.}
\label{fig:otl_stderr}
\end{figure}

Figure \ref{fig:otl_stderr} shows the (bootstrap) standard errors for the Monte Carlo estimates. Again, the bootstrap standard errors are generally lower than the relative errors from Figure \ref{fig:piston_err}. And all standard errors decay at the expected rate of $M^{-1/2}$. 

\subsection{Comparing rankings across metrics}

\begin{figure}[!h]
\centering
\subfloat[Piston model]{%
\includegraphics[width=0.45\textwidth]{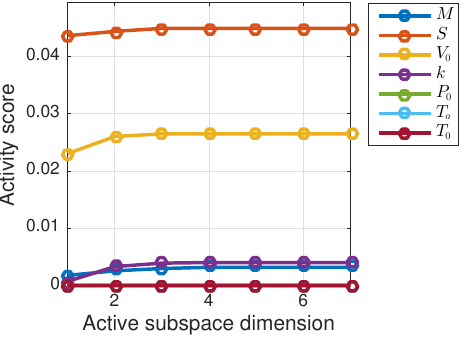}
\label{fig:alphan_piston}
}
\subfloat[Circuit model]{%
\includegraphics[width=0.45\textwidth]{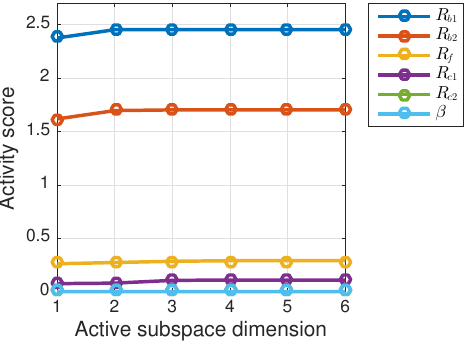}
\label{fig:alphan_otlcircuit}
}
\caption{Activity scores as a function of the active subspace dimension $n$. Recall that when $n=m$, the activity scores are the same as the derivative-based global sensitivity measures, as shown in Theorem \ref{thm:ascompare1}.}
\label{fig:activityscores}
\end{figure}

\noindent Figure \ref{fig:activityscores} shows the activity scores as a function of the active subspace dimension $n$ for each model. The rankings are consistent as $n$ increases with one exception: in the piston model, the rankings of parameters $M$ and $k$ swap as the dimension $n$ increases from 1 to 2. The activity scores converge quickly, which is due to the spectral decay of the eigenvalues weighting the squared eigenvector components in \eqref{eq:actscore}. Note that the lines representing the activity scores for parameters $P_0$ and $T_a$ in the piston model are behind the line for $T_0$, and all are near 0. Similarly, the lines representing the activity scores for parameter $R_{c2}$ are behind the line for $\beta$ in the circuit model.

\begin{figure}[!h]
\centering
\subfloat[Piston model]{%
\includegraphics[width=0.45\textwidth]{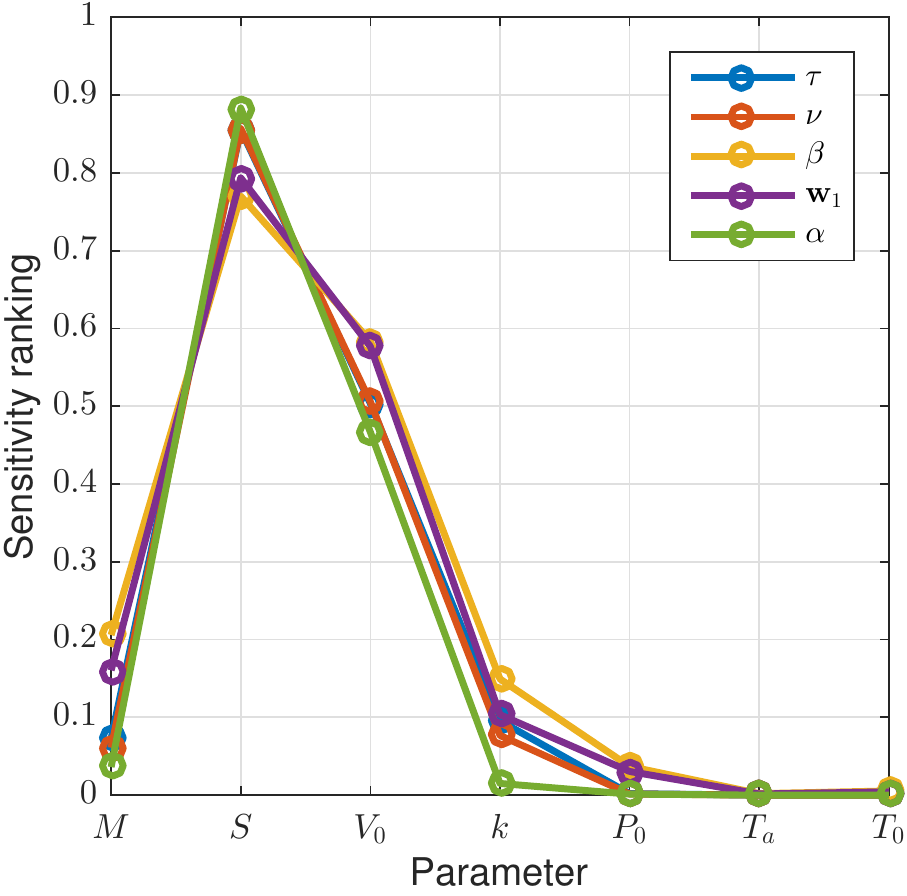}
\label{fig:alphan_piston}
}
\subfloat[Circuit model]{%
\includegraphics[width=0.45\textwidth]{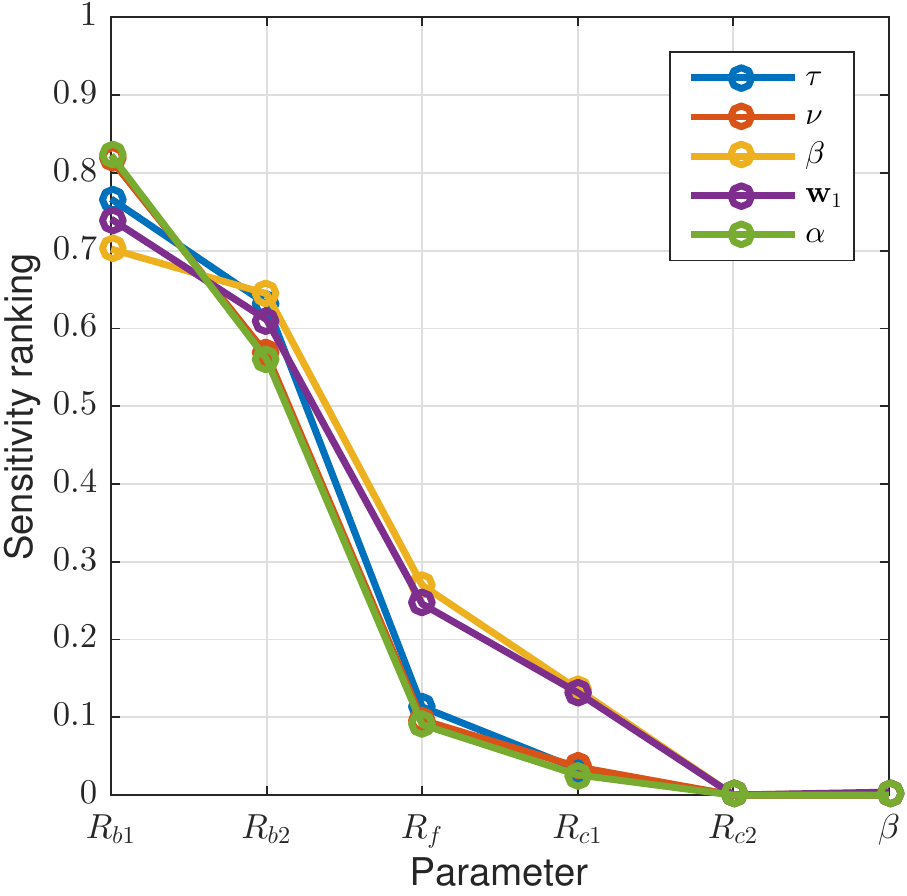}
\label{fig:alphan_otlcircuit}
}
\caption{Normalized sensitivity metrics allow visual comparison of importance rankings from the various sensitivity metrics for the two models. All metrics rank the parameters consistently for these two models.}
\label{fig:rankings}
\end{figure}

To visually compare rankings, we plot normalized versions of the reference values of the sensitivity metrics for each model. The normalization is as follows. For a sensitivity metric $\gamma_i$, we plot
\begin{equation}
\frac{|\gamma_i|}
{
\left(\sum_{i=1}^m \gamma_i^2\right)^{1/2}
},
\end{equation}
which (i) removes the signs from the linear model coefficients and the first eigenvector's components and (ii) scales the vector of numerical values to have norm 1. Figure \ref{fig:rankings} shows the normalized versions of the sensitivity metrics for both models. The parameters' rankings are consistent across the metrics for these two models. This supports our claim that the active subspace-based sensitivity metrics can be used like standard sensitivity metrics to measure the parameters' importance---at least, for these two models. And this complements the analysis from Section \ref{sec:compare}.

\section{Summary and remarks}
\label{sec:conc}

\noindent We propose two global sensitivity metrics derived from the eigenvectors and eigenvalues computed while testing an engineering model for an active subspace using its gradients. The two metrics are (i) the components of the first eigenvector and (ii) the \emph{activity scores}, which are the squared eigenvector components, linearly combined with the eigenvalues. We mathematically relate the activity scores to the Sobol' total sensitivity indices and the derivative-based global sensitivity measures. We discuss Gauss quadrature methods and Monte Carlo methods to estimate the new sensitivity metrics. And we demonstrate the new metrics on two algebraic engineering models. The models are simple enough to compute accurate quadrature-based reference values and study the behavior of the Monte Carlo estimates and their (bootstrap) standard errors estimates. For the two test models, all metrics identify the same important and unimportant variables, which supports our claim that the active subspace-based metrics behave similarly to the standard sensitivity metrics. While it is possible to construct models where the metrics rank parameters differently, since the metrics measure different characteristics of the model, we expect that this consistency will occur across many models in engineering practice. 

Our thesis is that the sensitivity metrics we derive from active subspaces behave similarly to the standard sensitivity metrics. We do not claim that the proposed metrics are computationally advantageous. In fact, the numerical experiments from Section \ref{sec:exp} suggest that the error in the Monte Carlo estimates of each metric are comparable for a fixed number $M$ of samples---particularly when comparing the activity scores to the derivative-based sensitivity measures. Similarly, the (bootstrap) standard errors are comparable across metrics, which suggests that the Monte Carlo estimates' variances are comparable for the same number of samples. Future work will analyze the error in the activity scores' Monte Carlo estimates as a function of the number of samples. Such analysis is beyond the scope of the current paper. 

Similar to the derivative-based global sensitivity measures, estimating the active subspace requires access to gradients. If gradients are not available and finite differences are not feasible, then one may seek to estimate active subspaces without gradients. Some initial work has been done along these lines~\cite{constantine2015sketching,Li2016}; we plan to explore how such gradient-free approaches lead to sensitivity metrics.

\section*{Acknowledgments}
This work is supported by the U.S. Department of Energy Office of Science, Office of Advanced Scientific Computing Research, Applied Mathematics program under Award Number DE-SC-0011077.





\bibliographystyle{elsarticle-num}
\bibliography{as-sensitivity-analysis}







\end{document}